\documentclass[psamsfonts]{amsart}
\usepackage[utf8]{inputenc}
\usepackage{amsfonts}
\usepackage{hyperref}

\usepackage{amsmath}
\usepackage{xcolor}
\usepackage{amsthm}
\usepackage{pdflscape}
\usepackage{mathrsfs}
\usepackage{bbold}
\usepackage{stmaryrd}
\usepackage[skip=5pt, indent=10pt]{parskip}

\usepackage{dsfont}
\usepackage{tikz}
\usepackage{amssymb}

\usepackage[all]{xy}

\usepackage{color}
\usepackage[shortlabels]{enumitem}
\usepackage{hyperref}
\usepackage{pgf,tikz}
\usepackage{mathrsfs}
\usetikzlibrary{arrows}
\usepackage{tikz-cd}
\usepackage{graphicx}
\usepackage{filecontents}

\newtheorem{thm}{Theorem}[section]
\newtheorem{cor}[thm]{Corollary}
\newtheorem{prop}[thm]{Proposition}
\newtheorem{lem}[thm]{Lemma}

\theoremstyle{definition}
\newtheorem{defn}[thm]{Definition}

\newtheorem{exmp}[thm]{Example}

\theoremstyle{remark}
\newtheorem{rem}[thm]{Remark}

\makeatletter
\let\c@equation\c@thm
\makeatother
\numberwithin{equation}{section}

\newcommand{\Z}{\mathbb{Z}}   
\newcommand{\Q}{\mathbb{Q}}

\newcommand{\X}{\mathbb{X}}

\newcommand{\wdo}{\dot{w}}

\newcommand{\extw}{\widetilde{W}}

\newcommand{\diag}{\text{diag}}
\newcommand{\val}{\textrm{val}} 

\newcommand{\ra}{\rightarrow}

\newcommand{\xra}{\xrightarrow}

\newcommand{\on}{\operatorname}

\newcommand{\tOm}{\widetilde{\Omega}}

\DeclareMathOperator{\id}{id}

\DeclareMathOperator{\rk}{rk}

\DeclareMathOperator{\Hom}{Hom}

\DeclareMathOperator{\Ad}{Ad}
\DeclareMathOperator{\aff}{aff}
\DeclareMathOperator{\Gal}{Gal}
\DeclareMathOperator{\Spec}{Spec}
\DeclareMathOperator{\chara}{char}

\DeclareMathOperator{\RHom}{\textbf{R} \kern -2pt \Hom}
\DeclareMathOperator{\Rhom}{\textbf{R}\kern -2pt \hom}
\DeclareMathOperator{\RuHom}{\textbf{R} \underline{\kern -0pt \Hom}}

\newcommand\cK{\mathcal{K}}
\newcommand\cL{\mathcal{L}}

\newcommand\cO{\mathcal{O}}
\newcommand\cP{\mathcal{P}}

\newcommand{\CC}{\mathbb{C}}

\newcommand{\RR}{\mathbb{R}}
\newcommand{\FF}{\mathbb{F}}

\newcommand\Ga{\mathbb{G}_a}

\newcommand{\der}{\textup{der}}

\numberwithin{equation}{section}

\newcommand{\fg}{\mathfrak{g}}

\renewcommand{\and}{\quad \on{and} \quad}

\newcommand\quash[1]{}

\newcommand{\lr}[1]{(\!(#1)\!)}

\renewcommand\c\circ

\renewcommand\a\alpha
\renewcommand\b\beta
\renewcommand\d\delta

\renewcommand\r\rho

\renewcommand{\t}{\tau}

\newcommand{\om}{\omega}
\newcommand{\Om}{\Omega}

\newcommand\Wa{W_{\aff}}

\newcommand\Nm{\on{Nm}}

\newcommand\pr{\textup{pr}}

\newcommand{\Id}{\mathrm{Id}}
\newcommand{\Daff}{\Delta_{\aff}}
\theoremstyle{plain}
\newtheorem{introthm}{Theorem}

\newtheorem{fact}[thm]{Fact}

\title{Sections of the Kottwitz homomorphism and affine pinnings}

\author{Yau Wing Li}
\address{School of Mathematics and Statistics, University of Melbourne, Parkville, VIC 3010, Australia}
\email{yauwingli57@gmail.com}
\subjclass[2020]{Primary 20G25; Secondary 20G07, 22E67}
\keywords{Kottwitz homomorphism, loop groups, affine pinnings, Tits lifts, Coxeter elements, centralizers of semisimple elements}

\begin{document}

\begin{abstract}
Let $G$ be a connected reductive group split over $\Z$. The Kottwitz homomorphism identifies the group of connected components of the loop group of $G$ with the group $\Om$ of length zero elements of the extended affine Weyl group. We show that this homomorphism always has a homomorphic section, and that the section can be chosen inside the normalizer of $T$ in $G(\Z[t,t^{-1}])$. We then ask whether the section can be chosen to preserve an affine pinning. This is possible if no simple factor of $G$ is of type $B_{2m}$ or $C_{2m}$. It is impossible for $(SL_2\times Sp_{2n})/\mu_2$ and $(SL_2\times Spin_{2n+1})/\mu_2$ when $n$ is even. We show that the obstructions for $G$ and for its Langlands dual group are dual to each other, and we compute them by a lattice computation with a Coxeter element. As an application, we show that the component group of the centralizer of any semisimple element lifts to the centralizer. Under the same condition on the simple factors, the lift can be chosen to preserve a pinning of the neutral component.
\end{abstract}

\maketitle

\section{Introduction}\label{sec:intro}

\subsection{Loop groups and their components}
Let $k$ be a field. For a linear algebraic group $G$ over $k$, the \emph{loop group} $LG$ is the functor $R\mapsto G(R(\!(t)\!))$ on $k$-algebras. It is a group ind-scheme with $k$-points $G(k(\!(t)\!))$. A reader who is not familiar with ind-schemes may simply think of the abstract group $G(k(\!(t)\!))$. Loop groups are the function-field analogues of $p$-adic groups. Their quotients by parahoric subgroups, such as the affine Grassmannian and the affine flag variety, are basic objects in geometric representation theory \cite{MV07,Be16,Zhu17}.

The name comes from topology. Over $\CC$, the group $G(\CC[t,t^{-1}])$ consists of algebraic maps $\CC^\times\ra G(\CC)$, and these restrict to loops in $G(\CC)$. The components of the space of loops in a connected Lie group $H$ are indexed by $\pi_1(H)$. Algebraic loop groups behave in the same way. Let $G$ be connected reductive and split, with split maximal torus $T$ and coroot lattice $Q^\vee\subset\X_*(T)$. Then the components of $LG$ are indexed by the algebraic fundamental group
\[
\pi_1(G)=\X_*(T)/Q^\vee
\]
\cite[Theorem~5.1]{PR08}. So $LG$ is disconnected unless $G$ is semisimple and simply connected.

Here are the basic examples. For $G=GL_n$, the map $g\mapsto\val(\det g)$ is a surjective homomorphism $GL_n(k(\!(t)\!))\ra\Z=\pi_1(GL_n)$. Its fibres are the components of $LGL_n$. For $SL_n$ the loop group is connected. For $PGL_n$ we have $PGL_n(k(\!(t)\!))=GL_n(k(\!(t)\!))/k(\!(t)\!)^\times$, and $\val\circ\det$ induces a surjection onto $\Z/n\Z=\pi_1(PGL_n)$.

In these examples, the projection onto the group of components has a homomorphic section. For $GL_n$ we can take $m\mapsto\diag(t^m,1,\dots,1)$. For $PGL_n$, consider
\[
\gamma=\begin{pmatrix}0&\Id_{n-1}\\ t&0\end{pmatrix}\in GL_n(k(\!(t)\!)).
\]
Then $\gamma^n=t\cdot\Id_n$ is a scalar matrix, and $\val(\det\gamma)=1$. So the image of $\gamma$ in $PGL_n(k(\!(t)\!))$ has order $n$, and $j\mapsto\gamma^j$ is a homomorphic section. The element $\gamma$ normalizes the standard Iwahori subgroup. It rotates the affine Dynkin diagram of type $A_{n-1}^{(1)}$, which is an $n$-cycle.

\subsection{The problem}
Does such a section exist for every split reductive group? A natural strategy is to lift each component to a \emph{monomial} element $\wdo\,t^{\mu}$. Here $\wdo$ lies in the normalizer $N_G(T)$, and $t^{\mu}=\mu(t)$ is a cocharacter evaluated at $t$. The matrix $\gamma$ is of this form.

The difficulty is that elements of the Weyl group $W=N_G(T)/T$ often cannot be lifted to $N_G(T)$ without increasing their order. For example, every element of $N_{SL_2}(T)$ over the non-trivial Weyl element has the form $\left(\begin{smallmatrix}0&a\\-a^{-1}&0\end{smallmatrix}\right)$, and its square is $-\Id$. Tits \cite{Ti66} constructed canonical lifts $n_w\in N_G(T)$ of all $w\in W$. They satisfy the braid relations, but $w\mapsto n_w$ is not a homomorphism.

For split almost simple groups over $p$-adic fields, Adrian constructed a homomorphic section from Tits' lifts \cite[Theorem~3.5]{Ad}. A new problem appears when several simple factors are glued together through the center. Consider
\[
G=(SL_2\times Sp_4)/\mu_2,
\]
where $\mu_2$ is embedded diagonally into the center. Then $\pi_1(G)\cong\Z/2\Z$. A monomial lift of the generator is the image in $G$ of a pair $(x,y)$, times a cocharacter evaluated at $t$. Here $x\in SL_2$ lies over the non-trivial Weyl element, and $y\in Sp_4$ lies over a certain involution $w$. The cocharacter cancels in the square, so the square of the lift is the image of $(x^2,y^2)$. Hence the lift has order $2$ if and only if
\[
y^2=x^2=-\Id.
\]
But Tits' lift satisfies $n_w^2=\Id$ (Example~\ref{exmp:sp4}). So $n_w$ must be replaced by another lift. Such a lift exists, but it does not respect the natural symmetries of the loop group (Theorem~\ref{thm:B}).

\subsection{Main results}
Let $G$ be a connected reductive group split over $\Z$, with split maximal torus $T$ and Weyl group $W$. Let $\extw=\X_*(T)\rtimes W$ be the extended affine Weyl group. We write its elements as $w\,t^{\mu}$. Let $\Om\subset\extw$ be the subgroup of elements of length zero. It is the stabilizer of the fundamental alcove, and it acts on the affine Dynkin diagram. The map $w\,t^{\mu}\mapsto[\mu]$ gives an isomorphism $\Om\cong\pi_1(G)$ \cite[Lemma~14]{HR08}.

For a field $k$, the projection of $LG_k$ onto its group of components is the Kottwitz homomorphism $\kappa_G$ \cite[\S7]{Ko97}, \cite[\S5]{PR08}. We obtain an exact sequence
\begin{equation}\label{eq:Kottwitz-seq}
1\ra LG_k^{\circ}\ra LG_k\xra{\kappa_G}\Om\ra 1,
\end{equation}
where $LG_k^{\circ}$ is the neutral component, and $\kappa_G(t^{\mu})=[\mu]$.

The monomial elements are defined over $\Z$. Every $n\in N_G(T)(\Z[t,t^{-1}])$ can be written as $n=\tau\,n_w\,t^{\mu}$ with $\tau\in T(\Z)$. The map $n\mapsto w\,t^{\mu}$ is a surjective homomorphism from $N_G(T)(\Z[t,t^{-1}])$ onto $\extw$, with kernel $T(\Z)\cong(\Z/2\Z)^{\rk T}$ (Lemma~\ref{lem:integral-ext}). Our first result says that this extension splits over $\Om$.

\begin{introthm}\label{thm:A}
Let $G$ be a connected reductive group split over $\Z$. There is a group homomorphism
\[
\psi\colon\Om\ra N_G(T)(\Z[t,t^{-1}])
\]
such that $\psi(\om)$ maps to $\om$ in $\extw$ for every $\om\in\Om$. Consequently, for every field $k$, the composite
\[
\Om\xra{\psi}N_G(T)(\Z[t,t^{-1}])\ra LG(k)
\]
is a homomorphic section of $\kappa_G$. In particular, $LG_k\cong LG_k^{\circ}\rtimes\Om$.
\end{introthm}

This is Theorem~\ref{thm:A-body}. Replacing $t$ by a uniformizer gives a homomorphic section of the Kottwitz homomorphism $G(F)\ra\Om$ for every non-archimedean local field $F$ (Corollary~\ref{cor:local-field}). This extends the results of Adrian from almost simple groups to all split reductive groups.

\subsection*{Pinned splittings}
A pinning of $G$ consists of a Borel subgroup $B\supset T$ and a non-zero root vector for each simple root. Automorphisms of $G$ that preserve a pinning correspond to automorphisms of the Dynkin diagram. The neutral component $LG^{\circ}$ has an analogous structure, an \emph{affine pinning}. It consists of an Iwahori subgroup $I\supset T$ and a root vector for each simple affine root (Definition~\ref{defn:affine-pinning}).

We say that $G$ is \emph{pinnable} over $k$ if \eqref{eq:Kottwitz-seq} has a splitting whose image preserves an affine pinning. Suppose $G$ is semisimple, and let $\tOm\subset LG(k)$ be the stabilizer of an affine pinning. Then there is an exact sequence (Proposition~\ref{prop:Z-extension} and Corollary~\ref{cor:surjective}; see also \cite[Proposition~3.14]{DLYZ})
\begin{equation}\label{eq:pinned-ext}
1\ra Z(G)(k)\ra\tOm\ra\Om\ra 1,
\end{equation}
and $G$ is pinnable if and only if \eqref{eq:pinned-ext} splits. When $G$ is adjoint, $Z(G)$ is trivial, so $G$ is pinnable.

\begin{introthm}\label{thm:B}
Let $G$ be a connected reductive group split over $\Z$, and let $k$ be a field.
\begin{enumerate}[(1)]
\item Suppose that $G$ is almost simple, or that no simple factor of $G$ is of type $B_{2m}$ or $C_{2m}$. Then $G$ is pinnable over $k$. Moreover, the section $\psi$ in Theorem~\ref{thm:A} can be chosen to preserve an affine pinning defined over $\Z$.
\item Suppose $\chara k\neq 2$, and let $n\geq 2$ be even. Let $G$ be $(SL_2\times Sp_{2n})/\mu_2$ or $(SL_2\times Spin_{2n+1})/\mu_2$, with $\mu_2$ embedded diagonally into the center. Then $\Om\cong Z(G)(k)\cong\Z/2\Z$, but every lift of the generator of $\Om$ that preserves an affine pinning has order $4$. In particular, $G$ is not pinnable over $k$.
\end{enumerate}
\end{introthm}

See Theorem~\ref{thm:B1-body} and Proposition~\ref{prop:B2-body}. So \eqref{eq:Kottwitz-seq} always splits, but a pinned splitting may not exist. The two groups in (2) are Langlands dual to each other. This is an instance of Theorem~\ref{thm:dual} below.

\subsection*{Duality}
Let $G$ be semisimple, and let $G^\vee$ be its Langlands dual group. Let $k$ be an algebraically closed field whose characteristic does not divide $|Z(G)|\cdot|\pi_1(G)|$. Then $Z(G)(k)$ is the Pontryagin dual of $\Om_{G^\vee}=\pi_1(G^\vee)$, and $\Om_G$ is the Pontryagin dual of $Z(G^\vee)(k)$. So the kernel and the quotient of \eqref{eq:pinned-ext} for $G$ are dual to the quotient and the kernel of \eqref{eq:pinned-ext} for $G^\vee$.

This says nothing about the extensions themselves, since an extension is not determined by its kernel and its quotient. For example, the groups $(SL_2\times Sp_4)/\mu_2$ and $(SL_2\times Sp_6)/\mu_2$ both have $Z(G)\cong\Om\cong\Z/2\Z$, but by Theorem~\ref{thm:B} the extension \eqref{eq:pinned-ext} splits only for the second one. The following theorem says that the extensions for $G$ and $G^\vee$ are dual as well.

\begin{introthm}\label{thm:dual}
Let $G$ be a split semisimple group, let $k$ be an algebraically closed field, and let $c\in W$ be a Coxeter element.
\begin{enumerate}[(1)]
\item Suppose that $\chara k\nmid|Z(G)|$. Then the extension \eqref{eq:pinned-ext} is isomorphic to
\[
0\ra Q^\vee/(1-c)\X_*(T)\ra\X_*(T)/(1-c)\X_*(T)\ra\X_*(T)/Q^\vee\ra0.
\]
\item Suppose that $\chara k\nmid|Z(G)|\cdot|\pi_1(G)|$. Then the extensions \eqref{eq:pinned-ext} for $G$ and for $G^\vee$ are Pontryagin dual to each other. In particular, $G$ is pinnable over $k$ if and only if $G^\vee$ is pinnable over $k$.
\end{enumerate}
\end{introthm}

Theorem~\ref{thm:dual} is proved in Section~\ref{sec:duality}, together with a more precise statement (Theorem~\ref{thm:duality}). The main idea is that the stabilizer of an affine pinning is the group of torsion points of an elliptic torus of Coxeter type, namely the centralizer of Kostant's cyclic element \cite{Kos59}. The author learned the use of Kostant's cyclic elements and Coxeter tori from the collaboration \cite{DLYZ} and discussions with Gurbir Dhillon, Zhiwei Yun and Xinwen Zhu. Part (1) turns pinnability into a lattice computation. For example, it shows that $(Sp_4\times Sp_6)/\mu_2$ is not pinnable (Example~\ref{exmp:sp4sp6}). This group is not covered by Theorem~\ref{thm:B}.

\subsection*{Centralizers}
Let $k$ be algebraically closed, and let $s\in G(k)$ be semisimple. Let $G_s$ be the centralizer of $s$, and let $\Om_s=\pi_0(G_s)$. The group $G_s$ is often disconnected. For example, in $PGL_2$ the centralizer of the image of $\diag(1,-1)$ is the normalizer of the diagonal torus. Such centralizers give rise to endoscopic groups.

\begin{introthm}\label{thm:cent}
Let $G$ be a connected reductive group split over $\Z$, let $k$ be an algebraically closed field, and let $s\in G(k)$ be semisimple.
\begin{enumerate}[(1)]
\item The extension
\[
1\ra G_s^{\circ}\ra G_s\ra\Om_s\ra 1
\]
splits. The splitting can be chosen to preserve a Borel subgroup of $G_s^{\circ}$ and a maximal torus in it.
\item If $G$ is pinnable over $k$, the splitting can be chosen to preserve a pinning of $G_s^{\circ}$. Then $\Om_s$ acts on $G_s^{\circ}$ by diagram automorphisms.
\item Suppose $\chara k\neq 2$, and let $G$ be one of the two groups in Theorem~\ref{thm:B}(2). Then there is an element $s\in G(k)$ of order $2$ such that no splitting preserves a pinning of $G_s^{\circ}$.
\end{enumerate}
\end{introthm}

See Theorem~\ref{thm:cent-body} and Proposition~\ref{prop:centralizer}. The splitting in (1) and (2) comes from the section of Theorem~\ref{thm:A} or Theorem~\ref{thm:B}. We restrict this section to the stabilizer of a point of the fundamental alcove, and then forget the translation parts: an element $\tau\,n_w\,t^{\mu}$ is sent to $\tau\,n_w$. For the section of Theorem~\ref{thm:A}, this simply means setting $t=1$. The key fact is that the walls of the alcove through this point give a basis of the root system of $G_s^{\circ}$ (Lemma~\ref{lem:walls}). Elements of infinite order are reduced to elements of finite order with the same centralizer (Lemma~\ref{lem:finite-order}).

\subsection*{The characteristic}
The results need the following assumptions on the field $k$.
\begin{itemize}
\item Theorem~\ref{thm:A}: none, since it is a statement over $\Z$.
\item Theorem~\ref{thm:B}: none for (1), and $\chara k\neq2$ for (2).
\item Theorem~\ref{thm:dual}: $k$ algebraically closed, with $\chara k\nmid|Z(G)|$ for (1) and $\chara k\nmid|Z(G)|\cdot|\pi_1(G)|$ for (2).
\item Theorem~\ref{thm:cent}: $k$ algebraically closed, and $\chara k\neq2$ for (3).
\end{itemize}
These assumptions cannot be dropped. In characteristic $2$, the groups in Theorem~\ref{thm:B}(2) are pinnable (Remark~\ref{rem:char2}). For Theorem~\ref{thm:dual}, let $\chara k=2$ and $G=SL_2$. Then the stabilizer of an affine pinning is $Z(G)(k)=1$, but the lattice extension in (1) has order $2$. Moreover, the stabilizers for $SL_2$ and for its dual group $PGL_2$ have orders $1$ and $2$, so (2) fails as well.

\subsection{Idea of the proof}
The main difficulty is to glue lifts across simple factors. Suppose $G=G^{sc}/K$ is semisimple, where $G^{sc}=\prod_iG_i$ is a product of simply connected almost simple groups. Then $\Om$ is a subgroup of $\prod_i\Om_{G_i^{ad}}$, and a lift of $\om\in\Om$ is built from lifts in the factors $G_i$. These lifts must be compatible modulo $K$. We achieve this by a uniform normalization. For $\om=w\,t^{\mu}\in\Om_{G_i^{ad}}$ of order $r$, we use a lift $\wdo\in N_{G_i}(T)(\Z)$ of $w$ with
\begin{equation}\label{eq:normalization}
\wdo^{\,r}=\exp(\mu)^{r(r-1)/2},
\end{equation}
where $\exp(\mu)$ is the central element of $G_i$ attached to $\om$.

The normalization is forced by type $A$: in $SL_r$, every lift of an $r$-cycle satisfies $\wdo^{\,r}=(-1)^{r-1}\Id=\exp(\mu)^{r(r-1)/2}$. It also has a conceptual meaning. Over $k\lr{t^{1/r}}$, the element $\wdo\,t^{\mu}$ lies in the simply connected group $G_i$, and \eqref{eq:normalization} says exactly that it is a Galois $1$-cocycle for $k\lr{t^{1/r}}/k\lr{t}$ (Lemma~\ref{lem:cocycle}). Products of cocycles are cocycles, and this is why normalized lifts can be glued across simple factors (Proposition~\ref{prop:semisimple} and Remark~\ref{rem:cocycle-gluing}).

Using Adrian's formula for powers of Tits' lifts \cite[Corollary~2.4]{Ad}, we show that Tits' lifts satisfy \eqref{eq:normalization} for the standard generators, except in types $B_{2m}$ and $C_{2m}$. In these types we correct Tits' lift by an element of $T(\Z)$. This correction is what breaks the pinning in Theorem~\ref{thm:B}(2). Finally, \eqref{eq:normalization} has to hold for every element of $\Om_{G_i^{ad}}$, not only for generators (Example~\ref{exmp:all-elements}).

\subsection{Motivation}
Our motivation comes from endoscopy for Hecke categories. Let $\cL$ be a character sheaf on $T$. It corresponds to a semisimple element $s$ of the dual torus $T^\vee\subset G^\vee$. Lusztig and Yun attached to $\cL$ an endoscopic group $H_{\cL}$, whose dual group is the neutral component of the centralizer of $s$ in $G^\vee$ \cite[\S1.2 and \S2.4]{LY20}. They showed that the neutral block of the monodromic Hecke category of $G$ is equivalent to the Hecke category of $H_{\cL}$ \cite[Theorem~1.3]{LY20}. They also treated the other blocks. These are indexed by the component group of the centralizer of $s$, and their monoidal structure involves a $3$-cocycle on this group \cite[Theorem~10.12]{LY20}, \cite{LY21}.

The affine analogue of the neutral block equivalence was proved in \cite[Theorem~1.1]{Li24}, and in a much more general setting in \cite{DLYZ}. In the affine case, the blocks are indexed by $\extw_{\cL}/\extw_{\cL}^{\circ}$ \cite[\S4.3]{Li24}. This group is an extension of the component group of the centralizer of $s$ by $\Om_{H_{\cL}}$. So the groups studied in this paper govern the non-neutral blocks. One possible approach to comparing the non-neutral blocks on the two sides of endoscopy is to lift these groups to the loop group, compatibly with pinnings. This paper studies when such lifts exist.

Lifting elements and subgroups of Weyl groups to torus normalizers is a classical theme \cite{Ti66,AH17,Lus18}. For connected reductive groups over non-archimedean local fields that split over an unramified extension, Ganapathy and He constructed a Tits group realizing an extension of the Iwahori--Weyl group by an elementary abelian $2$-group \cite{GH23}.

\subsection{Organization}
Section~\ref{sec:prelim} collects background on Tits' lifts, loop groups, fractional cocharacters and affine pinnings. In Section~\ref{sec:normalized} we construct normalized lifts in every simple type. Tits' lifts are normalized except in types $B_{2m}$ and $C_{2m}$, where we correct them (Table~\ref{tab:lifts}). In Section~\ref{sec:thmA} we show that normalized lifts glue along the center, and we prove Theorem~\ref{thm:A}; Remark~\ref{rem:cocycle-gluing} explains the gluing in terms of Galois cocycles. Section~\ref{sec:pinned} studies pinned splittings and proves Theorem~\ref{thm:B}. Section~\ref{sec:duality} proves the duality of the pinned extensions (Theorem~\ref{thm:dual}) and gives a lattice criterion for pinnability. Section~\ref{sec:centralizers} applies the splitting to centralizers of semisimple elements (Theorem~\ref{thm:cent}).

\subsection*{Acknowledgements}
I thank Dougal Davis and Kari Vilonen for helpful discussions, and Gurbir Dhillon, Zhiwei Yun and Xinwen Zhu for our collaboration and many valuable discussions. Most of this work was carried out at the University of Melbourne, with support from the University and the Australian Research Council grant FL200100141.

\section{Preliminaries}\label{sec:prelim}

In this section we fix notation and recall the facts we need. Sections~\ref{subsec:groups}--\ref{subsec:loop} and~\ref{subsec:W} are standard. Section~\ref{subsec:galois} introduces cocharacters with fractional exponents and the Galois action on them. This is where the normalization of Section~\ref{sec:normalized} comes from, and it is used again in Section~\ref{sec:duality}. Section~\ref{subsec:affine-pinning} recalls affine pinnings and their stabilizers.

\subsection{Split reductive groups}\label{subsec:groups}
Let $G$ be a connected reductive group split over $\Z$. We fix a split maximal torus $T$ and a Borel subgroup $B\supset T$. Let $\Phi$ be the set of roots, $\Phi^+$ the set of positive roots, and $\Delta$ the set of simple roots. For $\alpha\in\Phi$, let $\alpha^\vee$ be the coroot and $s_\alpha\in W$ the reflection. Let $Q^\vee\subset\X_*(T)$ be the coroot lattice. For an almost simple group, we label the simple roots as in the tables at the end of Chapter~VI of Bourbaki \cite[Plates~I--IX]{Bou02}, and we write $\theta$ for the highest root. We write $\varpi_i$ and $\varpi_i^\vee$ for the fundamental weights and coweights.

Split reductive groups over $\Z$ are classified by their root data \cite[Exp.~XXV]{SGA3}, and we use the standard models. For example, $Sp_{2n}$ is the group of automorphisms of $\Z^{2n}$ that preserve the standard alternating form, and $Spin_m$ is the spin group of the split quadratic form on $\Z^m$ \cite[Appendix~C]{Con14}. Both are simply connected. Recall that $Spin_5\cong Sp_4$, so the two groups in Theorem~\ref{thm:B}(2) coincide when $n=2$.

Suppose $G$ is semisimple. Let $P^\vee$ be the coweight lattice, so that
\[
Q^\vee\subset\X_*(T)\subset P^\vee.
\]
Let $\pi\colon G^{sc}\ra G$ be the simply connected cover, and put $K=\ker\pi$ and $T^{sc}=\pi^{-1}(T)$. Then $\X_*(T^{sc})=Q^\vee$, and $K$ is a finite subgroup scheme of the center $Z(G^{sc})\subset T^{sc}$. We write $G^{ad}$ for the adjoint group and $\Om^{ad}$ for its group $\Om$. So $\Om^{ad}\cong P^\vee/Q^\vee$.

For $x\in\X_*(T)\otimes\CC$, let $\exp(x)\in T(\CC)$ be its image under the exponential map, normalized so that $\exp(\lambda\otimes z)=\lambda(e^{2\pi iz})$. This is a surjective homomorphism with kernel $\X_*(T)$.

\begin{lem}\label{lem:exp}
Suppose $G$ is semisimple.
\begin{enumerate}[(1)]
\item The map $\mu\mapsto\exp(\mu)$ induces an isomorphism $P^\vee/Q^\vee\cong Z(G^{sc})(\CC)$.
\item This isomorphism identifies $\X_*(T)/Q^\vee$ with $K(\CC)$.
\end{enumerate}
\end{lem}

\begin{proof}
Let $x\in\X_*(T^{sc})\otimes\CC$. For $\alpha\in\Phi$ we have $\alpha(\exp x)=e^{2\pi i\langle\alpha,x\rangle}$. So $\exp(x)$ is central if and only if $\langle\alpha,x\rangle\in\Z$ for all $\alpha\in\Delta$, that is, if and only if $x\in P^\vee$. Since $\exp$ is surjective with kernel $Q^\vee$, this proves (1).

For (2), note that $K(\CC)$ is the set of $z\in Z(G^{sc})(\CC)$ with $\lambda(z)=1$ for all $\lambda\in\X^*(T)$. For $\mu\in P^\vee$ we have $\lambda(\exp\mu)=e^{2\pi i\langle\lambda,\mu\rangle}$. This equals $1$ for all $\lambda\in\X^*(T)$ if and only if $\mu\in\X_*(T)$.
\end{proof}

We will use two simple facts. First, every $z\in T(\CC)$ with $z^2=1$ is of the form $\lambda(-1)$ with $\lambda\in\X_*(T)$, so $z$ lies in $T(\Z)$. Second, let $\mu\in P^\vee$ and $r\geq1$ with $r\mu\in Q^\vee$. Then
\begin{equation}\label{eq:exp-power}
\exp(\mu)^{r(r-1)/2}=\big(r(r-1)\mu\big)(-1).
\end{equation}
Indeed, both sides equal $\exp\big(\tfrac{r-1}{2}\cdot r\mu\big)$. In particular, this element has order at most $2$. If $G$ is simply connected and $\lambda\in Q^\vee$, then $\lambda(-1)=1$ if and only if $\lambda\in 2Q^\vee$, if and only if $\langle\varpi_i,\lambda\rangle$ is even for all $i$.

\subsection{Pinnings and Tits' lifts}\label{subsec:tits}
Tits' lifts are our starting point for lifting Weyl group elements. We recall their definition and two formulas.

We fix a pinning, that is, isomorphisms $u_\alpha\colon\Ga\ra U_\alpha$ for $\alpha\in\Delta$, defined over $\Z$. They extend to a realization $(u_\alpha)_{\alpha\in\Phi}$ of $\Phi$ in $G$ \cite[\S8.1]{Sp}. For $\alpha\in\Phi$ we put
\[
n_\alpha=u_\alpha(1)\,u_{-\alpha}(-1)\,u_\alpha(1)\in N_G(T)(\Z).
\]
It lifts $s_\alpha$, and $n_\alpha^2=\alpha^\vee(-1)$. For $w\in W$, choose a reduced expression $w=s_{\alpha_{i_1}}\cdots s_{\alpha_{i_m}}$ and put $n_w=n_{\alpha_{i_1}}\cdots n_{\alpha_{i_m}}$. This does not depend on the reduced expression \cite{Ti66}, \cite[\S9.3]{Sp}. We call $n_w$ the \emph{Tits lift} of $w$.

\begin{lem}[{\cite[Proposition~9.3.5]{Sp}}]\label{lem:springer}
Let $\alpha,\beta\in\Delta$ and $w\in W$ with $w(\alpha)=\beta$. Then $n_w\,u_\alpha(c)\,n_w^{-1}=u_\beta(c)$.
\end{lem}

\begin{lem}[{\cite[Proposition~2.2 and Corollary~2.4]{Ad}}]\label{lem:adrian}
\begin{enumerate}[(1)]
\item For $u,v\in W$, let $F(u,v)=\{\alpha\in\Phi^+: v(\alpha)\in-\Phi^+,\ uv(\alpha)\in\Phi^+\}$. Then
\[
n_u n_v=n_{uv}\prod_{\alpha\in F(u,v)}\alpha^\vee(-1).
\]
\item For $w\in W$ and $i\geq1$, let $I_i(w)=F(w,w^i)$. Then for $m\geq1$,
\[
n_w^m=n_{w^m}\prod_{i=1}^{m-1}\prod_{\alpha\in I_i(w)}\alpha^\vee(-1).
\]
\end{enumerate}
\end{lem}

If $w^r=1$, it follows that $n_w^r=\gamma_w(-1)$, where
\[
\gamma_w=\sum_{i=1}^{r-1}\sum_{\alpha\in I_i(w)}\alpha^\vee\in Q^\vee.
\]

\begin{exmp}\label{exmp:sp4}
Let $V$ be a symplectic space with basis $v_1,v_2,v_3,v_4$ such that $\langle v_1,v_4\rangle=\langle v_2,v_3\rangle=1$ and $\langle v_i,v_j\rangle=0$ if $i+j\neq5$. Let $T\subset Sp(V)$ be the torus of diagonal matrices $\diag(b,c,c^{-1},b^{-1})$, and let $L_1,L_2$ be the characters $b$ and $c$. The simple roots are $\alpha_1=L_1-L_2$ (short) and $\alpha_2=2L_2$ (long). For a suitable pinning,
\[
n_{\alpha_1}=\left(\begin{array}{rrrr}
0&1&0&0\\ -1&0&0&0\\ 0&0&0&-1\\ 0&0&1&0
\end{array}\right),\qquad
n_{\alpha_2}=\left(\begin{array}{rrrr}
1&0&0&0\\ 0&0&1&0\\ 0&-1&0&0\\ 0&0&0&1
\end{array}\right).
\]
Let $w=s_{\alpha_2}s_{\alpha_1}s_{\alpha_2}=s_{L_1+L_2}$. Then
\[
n_w=n_{\alpha_2}n_{\alpha_1}n_{\alpha_2}=\left(\begin{array}{rrrr}
0&0&1&0\\ 0&0&0&-1\\ 1&0&0&0\\ 0&-1&0&0
\end{array}\right),
\]
and $n_w^2=\Id$. This is the involution $w$ of the introduction.
\end{exmp}

\subsection{Loop groups and the Kottwitz homomorphism}\label{subsec:loop}
We now recall the loop group and its components.

Let $k$ be a field. Put $\cK=k(\!(t)\!)$ and $\cO=k[[t]]$. Let $LG$ and $L^+G$ be the loop group and the positive loop group over $k$, so
\[
LG(R)=G(R(\!(t)\!)),\qquad L^+G(R)=G(R[[t]]).
\]
Let $I\subset L^+G$ be the standard Iwahori subgroup, the preimage of $B$ under $t\mapsto0$.

\begin{thm}[{\cite[Theorem~5.1, Corollary~5.3]{PR08}}]\label{thm:PR}
Suppose $k$ is algebraically closed. Then the Kottwitz homomorphism induces a bijection $\pi_0(LG)\cong\pi_1(G)$.
\end{thm}

For an arbitrary field $k$, each component of $LG_{\bar k}$ contains a point $t^{\mu}$ with $\mu\in\X_*(T)$. Hence it is defined over $k$, and $\pi_0(LG_k)\cong\pi_1(G)$. We write
\[
\kappa_G\colon LG(k)\ra\pi_1(G)\cong\Om
\]
for the resulting homomorphism. It is normalized by $\kappa_G(t^{\mu})=[\mu]$. Its kernel is the group of $k$-points of the neutral component $LG_k^{\circ}$.

We only use loop groups over fields. Over $\Z$, the structure of $LG$ is more delicate. For example, if $p$ divides $|\pi_1(G_{\der})|$, then $LG_{\FF_p}$ is not reduced \cite[Remark~6.4]{PR08}, \cite{HLR24}.

The monomial elements live in the following group.

\begin{lem}\label{lem:integral-ext}
\begin{enumerate}[(1)]
\item The group $T(\Z[t,t^{-1}])$ is equal to $T(\Z)\times t^{\X_*(T)}$. Moreover, $T(\Z)\cong(\Z/2\Z)^{\rk T}$.
\item Let $n\in N_G(T)(\Z[t,t^{-1}])$. Then $n$ can be written uniquely as $n=\tau\,n_w\,t^{\mu}$ with $\tau\in T(\Z)$, $w\in W$ and $\mu\in\X_*(T)$. The map $n\mapsto w\,t^{\mu}$ is a surjective homomorphism from $N_G(T)(\Z[t,t^{-1}])$ onto $\extw$, with kernel $T(\Z)$.
\end{enumerate}
\end{lem}

\begin{proof}
(1) Since $T$ is split, $T(R)=\Hom(\X^*(T),R^\times)$. The units of $\Z[t,t^{-1}]$ are $\pm t^m$.

(2) The quotient $N_G(T)/T$ is the constant group $W$, and $\Spec\Z[t,t^{-1}]$ is connected. So $n$ has a well-defined image $w\in W$, and $n\,n_w^{-1}\in T(\Z[t,t^{-1}])$. By (1), $n=\tau\,t^{\nu}n_w=\tau\,n_w\,t^{w^{-1}\nu}$ for some $\tau\in T(\Z)$ and $\nu\in\X_*(T)$. Uniqueness is clear. For the homomorphism property, we use the relation $t^{\mu}n_{v}=n_{v}\,t^{v^{-1}\mu}$. We also use that $n_wn_v\,n_{wv}^{-1}$ lies in $T(\Z)$, by Lemma~\ref{lem:adrian}(1).
\end{proof}

\subsection{Fractional cocharacters and the Galois action}\label{subsec:galois}
The elements of $\Om^{ad}$ have translation parts in the coweight lattice $P^\vee$, which is in general larger than $\X_*(T)$. To lift them to a simply connected group, we need cocharacters with fractional exponents. These are defined over a ramified extension of $\cK$, and the Galois group of this extension enters the picture.

Let $k$ be algebraically closed. Fix a compatible system of primitive roots of unity $\zeta_r\in k$, for $r$ prime to $\chara k$, so that $\zeta_{rs}^{\,s}=\zeta_r$. For $k=\CC$ we take $\zeta_r=e^{2\pi i/r}$. For such $r$, put
\[
\cK_r=k\lr{t^{1/r}}.
\]
This is a Galois extension of $\cK$ with cyclic Galois group $\Gamma_r=\Gal(\cK_r/\cK)$. It is generated by the element $\gamma$ with $\gamma(t^{1/r})=\zeta_r\,t^{1/r}$.

If $s$ divides $r$, then $\cK_s\subset\cK_r$ via $t^{1/s}=(t^{1/r})^{r/s}$. Since $\zeta_r^{\,r/s}=\zeta_s$, we get
\[
\gamma(t^{1/s})=\zeta_r^{\,r/s}\,t^{1/s}=\zeta_s\,t^{1/s}.
\]
So $\gamma$ acts on $\cK_s$ as the generator of $\Gamma_s$. In particular, $\gamma^s$ acts trivially on $\cK_s$.

Let $S$ be a split torus, and let $\mu\in\X_*(S)\otimes\Q$ with $r\mu\in\X_*(S)$. We put
\[
t^\mu=(r\mu)(t^{1/r})\in S(\cK_r).
\]
This does not depend on $r$. Its Galois conjugate is
\begin{equation}\label{eq:monodromy}
\gamma(t^\mu)=\exp(\mu)\,t^\mu,\qquad\text{where }\exp(\mu):=(r\mu)(\zeta_r)\in S(k).
\end{equation}
We call $\exp(\mu)$ the \emph{monodromy} of $t^\mu$. For $k=\CC$, this agrees with the exponential map of Section~\ref{subsec:groups}.

Now let $G$ be semisimple, and apply this to $S=T^{sc}$ and $\mu\in P^\vee$. Then $\exp(\mu)$ is central in $G^{sc}$, because $\alpha(\exp\mu)=\zeta_r^{\,r\langle\alpha,\mu\rangle}=1$ for every root $\alpha$. Lemma~\ref{lem:exp} holds over $k$ with the same proof: the map $\exp$ induces $P^\vee/Q^\vee\cong Z(G^{sc})(k)$, and $\X_*(T)/Q^\vee$ corresponds to $K(k)$. The element $\exp(\mu)^{r(r-1)/2}$ of \eqref{eq:exp-power} does not depend on the choice of $\zeta_r$.

Conjugation by $t^\mu$ preserves $LG^{sc}=G^{sc}(\cK)$, because
\[
t^\mu u_\alpha(f)\,t^{-\mu}=u_\alpha\big(t^{\langle\alpha,\mu\rangle}f\big)\quad\text{and}\quad\langle\alpha,\mu\rangle\in\Z.
\]
If $\dot x\in N_{G^{sc}}(T^{sc})$ lifts $w\in W$, then $t^\mu\dot x=\dot x\,t^{w^{-1}\mu}$.

The monodromy also describes the Kottwitz homomorphism of a semisimple group. In other words, $\kappa_G$ is a Kummer map.

\begin{lem}\label{lem:kummer}
Let $G=G^{sc}/K$ be semisimple, and let $k$ be algebraically closed with $\chara k\nmid|K|$. Let $r$ be a multiple of the exponent of $\Om$ that is prime to $\chara k$.
\begin{enumerate}[(1)]
\item Every $y\in G(\cK)$ lifts to an element $\tilde y\in G^{sc}(\cK_r)$.
\item For such a lift, $\tilde y^{-1}\gamma(\tilde y)$ lies in $K(k)$, and it does not depend on the choice of $\tilde y$.
\item Under the isomorphism $K(k)\cong\X_*(T)/Q^\vee\cong\Om$, the element $\tilde y^{-1}\gamma(\tilde y)$ equals $\kappa_G(y)$.
\end{enumerate}
\end{lem}

\begin{proof}
We first show that $G(\cK)$ is generated by $\pi(G^{sc}(\cK))$ and the elements $\mu(t)$ with $\mu\in\X_*(T)$. The group $LG^\circ(k)$ is generated by $T(\cO)$ and the affine root subgroups. These root subgroups lie in $\pi(G^{sc}(\cK))$. Moreover $T(\cO)\subset\pi(T^{sc}(\cO))$, because $|K|$ is prime to $\chara k$. So $LG^\circ(k)\subset\pi(G^{sc}(\cK))$. Since every component of $LG$ contains some $\mu(t)$, the claim follows.

Now $\mu(t)$ lifts to $t^\mu\in T^{sc}(\cK_r)$, since $r\mu\in Q^\vee$. With the claim, this proves (1). For (2), the element $\pi(\tilde y)=y$ is Galois invariant, so $\tilde y^{-1}\gamma(\tilde y)$ lies in $K(\cK_r)=K(k)$. Two lifts differ by an element of $K(k)$, which is central and Galois invariant. So the element does not depend on the lift.

For (3), the map $\delta\colon y\mapsto\tilde y^{-1}\gamma(\tilde y)$ is a homomorphism $G(\cK)\ra K(k)$, because $K$ is central. It vanishes on $\pi(G^{sc}(\cK))$. By \eqref{eq:monodromy}, $\delta(\mu(t))=\exp(\mu)$, which corresponds to $[\mu]=\kappa_G(\mu(t))$. Since $\delta$ and $\kappa_G$ agree on generators, they are equal.
\end{proof}

\subsection{The extended affine Weyl group}\label{subsec:W}
We recall the action of $\extw$ on the apartment, the group $\Om$, and its standard generators.

Let $V=\X_*(T)\otimes\RR$. The group $\extw$ acts on $V$ by affine transformations,
\[
w\,t^{\mu}\cdot x=w(x-\mu).
\]
An \emph{affine root} is a pair $a=\alpha+m\delta$ with $\alpha\in\Phi$ and $m\in\Z$, viewed as the affine function $a(x)=\langle\alpha,x\rangle+m$ on $V$. We write $\pi(a)=\alpha$ for its gradient. The group $\extw$ acts on affine roots by $(g\cdot a)(x)=a(g^{-1}x)$. Explicitly,
\[
w\,t^{\mu}\cdot(\alpha+m\delta)=w\alpha+(m+\langle\alpha,\mu\rangle)\delta.
\]
The root subgroup of $a$ is $U_a=\{u_\alpha(c\,t^m):c\in k\}\subset LG$. If $n\in N_G(T)(\cK)$ maps to $g\in\extw$, then $nU_an^{-1}=U_{g\cdot a}$. This is why we use the above action.

Let
\[
A=\{x\in V:0<\langle\alpha,x\rangle<1\text{ for all }\alpha\in\Phi^+\}
\]
be the fundamental alcove. The affine roots that are positive on $A$ are the $\alpha+m\delta$ with $m\geq0$ if $\alpha>0$, and $m\geq1$ if $\alpha<0$. The set $\Daff$ of simple affine roots consists of the simple roots and the roots $\delta-\theta$, one for each simple factor. The Iwahori subgroup $I$ is generated by $T(\cO)$ and the $U_a$ with $a$ positive on $A$.

The affine Weyl group $\Wa=Q^\vee\rtimes W$ is generated by the reflections in the walls of $A$. It acts simply transitively on the set of alcoves. Let $\Om$ be the stabilizer of $A$ in $\extw$. Then
\[
\extw=\Wa\rtimes\Om,
\]
the group $\Om$ permutes $\Daff$, and $w\,t^{\mu}\mapsto[\mu]$ is an isomorphism $\Om\cong\pi_1(G)$ \cite{IM65}, \cite[Lemma~14]{HR08}. Every element of $N_G(T)(\cK)$ that maps to $\Om$ normalizes $I$.

Let $\pr\colon\extw\ra W$ be the projection $w\,t^{\mu}\mapsto w$. We call $\pr(g)$ the \emph{finite part} of $g\in\extw$. If $t^{\mu}\in\Om$, then $\langle\alpha,\mu\rangle=0$ for all $\alpha$. So the kernel of $\pr|_{\Om}$ is $\X_*(Z(G)^{\circ})$, and $\pr|_{\Om}$ is injective when $G$ is semisimple.

Suppose $G$ is almost simple. A coweight $\varpi_j^\vee$ is minuscule if the coefficient of $\alpha_j$ in $\theta$ is $1$. For such $j$, let $w_{0,j}$ be the longest element of the subgroup of $W$ generated by the $s_{\alpha_i}$ with $i\neq j$, and let $w_0$ be the longest element of $W$. By \cite[Proposition~1.18]{IM65}, the non-trivial elements of $\Om^{ad}$ are
\[
\om_j=w_0w_{0,j}\,t^{\varpi_j^\vee},\qquad \varpi_j^\vee\text{ minuscule}.
\]
In particular, $\exp(\varpi_j^\vee)\in Z(G^{sc})$ is the element attached to $\om_j$ by Lemma~\ref{lem:exp}.

\subsection{Affine pinnings}\label{subsec:affine-pinning}
Finally, we recall affine pinnings and their stabilizers. These are the objects of Theorem~\ref{thm:B} and Theorem~\ref{thm:dual}.

\begin{defn}\label{defn:affine-pinning}
An \emph{affine pinning} of $LG_k^{\circ}$ is a tuple $\cP=(I,T,(x_a)_{a\in\Daff})$, where $x_a\in U_a(k)$ is non-trivial for every $a\in\Daff$. Its \emph{stabilizer} $\tOm_{\cP}$ is the group of $g\in LG(k)$ that normalize $I$ and $T$ and satisfy
\[
g\,x_a\,g^{-1}\in\{x_b:b\in\Daff\}\qquad\text{for all }a.
\]
\end{defn}

For example, a pinning of $G$ and a choice of signs give the \emph{standard} affine pinnings, with $x_\alpha=u_\alpha(1)$ for $\alpha\in\Delta$ and $x_{\delta-\theta}=u_{-\theta}(\pm t)$.

\begin{defn}\label{defn:pinnable}
We say that $G$ is \emph{pinnable over $k$} if there are an affine pinning $\cP$ and a homomorphic section $\sigma$ of $\kappa_G$ with $\sigma(\Om)\subset\tOm_{\cP}$.
\end{defn}

\begin{prop}\label{prop:Z-extension}
Let $G$ be semisimple, and let $\cP$ be an affine pinning of $LG_k$. Then $\tOm_\cP$ is a finite subgroup of $LG(k)$, and the restriction of $\kappa_G$ to $\tOm_\cP$ has kernel $Z(G)(k)$. In particular, if $\kappa_G(\tOm_\cP)=\Om$, then
\[
1\ra Z(G)(k)\ra\tOm_\cP\xra{\ \kappa_G\ }\Om\ra1
\]
is a central extension of finite groups. Its splittings are exactly the homomorphic sections of $\kappa_G$ with image in $\tOm_\cP$.
\end{prop}

\begin{proof}
Clearly $Z(G)(k)$ lies in the kernel. Conversely, let $g\in\tOm_{\cP}$ with $\kappa_G(g)=0$. Then $g$ normalizes $T$, so it normalizes the centralizer $LT$ of $T$. Since $g$ also normalizes $I$, its image in $\extw$ lies in $\Om$. As $\kappa_G(g)=0$, this image is trivial, and $g\in T(\cO)$. Now $g$ fixes every $x_a$. So $\alpha(g)=1$ for $\alpha\in\Delta$, and $\theta(g)=1$ for each highest root. Hence $g\in Z(G)(\cO)$. Since $Z(G)$ is a finite multiplicative group scheme and $\cO$ is a henselian local domain, we get $Z(G)(\cO)=Z(G)(k)$. Finally, $\tOm_\cP$ is finite because $Z(G)(k)$ and $\Om$ are finite.
\end{proof}

The condition $\kappa_G(\tOm_\cP)=\Om$ is fulfilled over every field $k$ for a suitable affine pinning $\cP$ defined over $\Z$, and for every affine pinning if $k$ is algebraically closed (Corollary~\ref{cor:surjective}). Moreover, if $\chara k\nmid|Z(G)|$, then $\tOm_\cP$ is commutative (Theorem~\ref{thm:duality}(1), applied over $\bar k$). For connected reductive $G$, see also \cite[Proposition~3.14]{DLYZ}, where $\tOm_\cP$ is denoted $M_\phi$.

To compare different affine pinnings, we use the following lemma. For an affine pinning $\cP=(I,T,(x_a))$, let $u_a\colon\Ga\ra U_a$ be the isomorphism with $u_a(1)=x_a$. For $c\in(k^\times)^{\Daff}$, put
\[
\cP_c=(I,T,(u_a(c_a))_{a\in\Daff}).
\]
For $a\in\Daff$ and $\tau\in T(k)$, we write $a(\tau)$ for $\pi(a)(\tau)$.

\begin{lem}\label{lem:conj-pinning}
Suppose $G$ is semisimple and $k$ is algebraically closed. Let $\cP$ be an affine pinning, and let $g\in\tOm_{\cP}$ with $\kappa_G(g)=\om$. Let $\cP'$ be another affine pinning, and let $g'\in\tOm_{\cP'}$ with $\kappa_G(g')=\om$. Then there are $h\in LG(k)$ and $c\in(k^\times)^{\Daff}$ such that
\[
h\cP'h^{-1}=\cP_c\quad\text{and}\quad hg'h^{-1}=g\tau\ \text{ with }\tau\in T(k).
\]
Moreover, $a(\tau)=c_{\om(a)}/c_a$ for all $a\in\Daff$.
\end{lem}

\begin{proof}
All Iwahori subgroups of $LG(k)$ are conjugate, and all maximal tori of $I$ are conjugate under $I$. So we may assume that $\cP'=\cP_c$ for some $c$. Then $g$ and $g'$ normalize $I$ and $T$, and they have the same image $\om$ in $\Om\subset\extw$. Hence $g'=g\tau$ with $\tau\in T(\cO)$.

Since $g$ normalizes $T$, conjugation by $g$ is linear in the coordinates $u_a$. As $g$ preserves $\cP$, it sends $u_a(x)$ to $u_{\om(a)}(x)$. So conjugation by $g'$ sends $u_a(c_a)$ to $u_{\om(a)}(a(\tau)c_a)$, and this equals $u_{\om(a)}(c_{\om(a)})$. Thus $a(\tau)=c_{\om(a)}/c_a$ is a constant for every $a\in\Daff$. Hence the image of $\tau$ in $T^{ad}(\cO)$ is constant. Since $Z(G)(\cO)=Z(G)(k)$, we get $\tau\in T(k)$.
\end{proof}

\section{Normalized lifts for almost simple groups}\label{sec:normalized}

In this section $G$ is almost simple and simply connected. Each element $\om=w\,t^\mu$ of $\Om^{ad}$ gives a Weyl group element $w$. We want lifts of $w$ to $N_G(T)(\Z)$ whose $r$-th powers are controlled, where $r$ is the order of $\om$.

The right condition is \eqref{eq:normalization} from the introduction. We first explain where it comes from: it says that a certain element over $\cK_r$ is a Galois $1$-cocycle (Lemma~\ref{lem:cocycle}), and type $A$ satisfies it automatically (Remark~\ref{rem:typeA}). Then we show that every type admits such lifts for the standard generators of $\Om^{ad}$ (Proposition~\ref{prop:generators}). The proof goes through the types one by one, using Adrian's formula for the powers of Tits' lifts. Finally, we pass from generators to all elements of $\Om^{ad}$ (Lemma~\ref{lem:powers}).

By Lemma~\ref{lem:exp}, the map $\om=w\,t^{\mu}\mapsto\exp(\mu)$ is an isomorphism $\Om^{ad}\cong Z(G)(\CC)$. Since $\pr$ is injective on $\Om^{ad}$, the elements $\om$ and $w$ have the same order.

\begin{defn}\label{defn:normalized}
Let $\om=w\,t^{\mu}\in\Om^{ad}$ have order $r$. A \emph{normalized lift} of $\om$ is an element $\dot x\in N_G(T)(\Z)$ lifting $w$ such that
\[
\dot x^{\,r}=\exp(\mu)^{r(r-1)/2}.
\]
\end{defn}

By \eqref{eq:exp-power}, the right hand side is a central element of order at most $2$. In particular, the condition does not depend on the choice of roots of unity in Section~\ref{subsec:galois}.

\subsection{The cocycle condition}
Here is the conceptual meaning of the normalization. Over $\cK_r$, the element $\dot x\,t^\mu$ lies in the simply connected group $G$, and normalized lifts are exactly those for which it is a Galois $1$-cocycle. In the next lemma, $k$ is algebraically closed, and $G$ may be any simply connected semisimple group.

\begin{lem}\label{lem:cocycle}
Let $\om=w\,t^\mu\in\Om^{ad}$, let $r$ be a multiple of its order with $\chara k\nmid r$, and let $\dot x\in N_G(T)(k)$ lift $w$. Put $g=\dot x\,t^\mu\in G(\cK_r)$. Then $\gamma(g)=\exp(\mu)\,g$, and
\[
N_r(g):=g\,\gamma(g)\cdots\gamma^{r-1}(g)=\exp(\mu)^{r(r-1)/2}\,\dot x^{\,r}.
\]
In particular, if $r$ is the order of $\om$, then $\dot x$ is a normalized lift if and only if $\gamma\mapsto g$ defines a $1$-cocycle of $\Gamma_r$ with values in $G(\cK_r)$.
\end{lem}

\begin{proof}
By \eqref{eq:monodromy}, $\gamma(g)=\dot x\exp(\mu)t^\mu=\exp(\mu)\,g$, since $\exp(\mu)$ is central. Hence $\gamma^j(g)=\exp(\mu)^jg$, and
\[
N_r(g)=\exp(\mu)^{0+1+\dots+(r-1)}\,g^r.
\]
Moreover, $g^r=\dot x^{\,r}\,t^{\nu}$, where $t^\nu$ is the translation part of $\om^r=1$. So $g^r=\dot x^{\,r}$, which gives the formula.

Recall that an element $g$ defines a $1$-cocycle of the cyclic group $\Gamma_r$ if and only if $N_r(g)=1$. Since $\exp(\mu)^{r(r-1)/2}$ has order at most $2$, this means
\[
\dot x^{\,r}=\exp(\mu)^{r(r-1)/2}.
\]
\end{proof}

\begin{rem}[type $A$]\label{rem:typeA}
In type $A$, every lift is a cocycle lift. Indeed, let $G=SL_r$ and let $\gamma_A=\left(\begin{smallmatrix}0&\Id_{r-1}\\t&0\end{smallmatrix}\right)$, which is Galois invariant with $\gamma_A^r=t\cdot\Id$. A lift of the generator of $\Om^{ad}$ is $g=\epsilon\,t^{-1/r}\gamma_A$, with $\epsilon^r=(-1)^{r-1}$ so that $\det g=1$. The scalar $t^{-1/r}$ is central, so
\[
N_r(g)=\epsilon^r\cdot\Nm_{\cK_r/\cK}(t^{-1/r})\cdot\gamma_A^r,\qquad \Nm_{\cK_r/\cK}(t^{1/r})=\zeta_r^{\,0+1+\dots+(r-1)}\,t=(-1)^{r-1}t.
\]
Hence $N_r(g)=1$. So the exponent $r(r-1)/2$ in Definition~\ref{defn:normalized} is the exponent of $\zeta_r$ in the norm of $t^{1/r}$.
\end{rem}

\subsection{Normalized lifts of generators}
We now find normalized lifts of the standard generators of $\Om^{ad}$. We first reduce the normalization to a congruence between coroots. This congruence can then be checked with Lemma~\ref{lem:adrian}.

\begin{lem}\label{lem:criterion}
Let $\om=w\,t^{\mu}\in\Om^{ad}$ have order $r$, and let $\lambda\in Q^\vee$. Then $n_w\,\lambda(-1)$ is a normalized lift of $\om$ if and only if
\[
\gamma_w+\sum_{j=0}^{r-1}w^j(\lambda)-r(r-1)\mu\in2Q^\vee.
\]
\end{lem}

\begin{proof}
Put $\tau=\lambda(-1)$. Since $n_w\tau n_w^{-1}=w(\tau)$, we have
\[
(n_w\tau)^r=\Big(\prod_{j=1}^rw^j(\tau)\Big)\,n_w^r.
\]
By Lemma~\ref{lem:adrian}, $n_w^r=\gamma_w(-1)$. Hence $(n_w\tau)^r=\big(\gamma_w+\sum_jw^j(\lambda)\big)(-1)$. By \eqref{eq:exp-power}, $\exp(\mu)^{r(r-1)/2}=(r(r-1)\mu)(-1)$. Finally, $\nu(-1)=1$ for $\nu\in Q^\vee$ if and only if $\nu\in2Q^\vee$.
\end{proof}

\begin{prop}\label{prop:generators}
Table~\ref{tab:lifts} lists generators $\om_j=w\,t^{\varpi_j^\vee}$ of $\Om^{ad}$, and normalized lifts of them. Here $\om_j$ is the element of Section~\ref{subsec:W} attached to the minuscule coweight $\varpi_j^\vee$, so $w=w_0w_{0,j}$; we use the Bourbaki labelling. In type $D_\ell$ with $\ell$ even, the two lifts commute. In types $E_8$, $F_4$ and $G_2$, the group $\Om^{ad}$ is trivial.
\end{prop}

\begin{table}[ht]
\renewcommand{\arraystretch}{1.2}
\begin{tabular}{lllll}
type & $G$ & generator & $r$ & normalized lift\\ \hline
$A_\ell$ & $SL_{\ell+1}$ & $\om_1$ & $\ell+1$ & $n_w$ (any lift)\\
$B_\ell$ & $Spin_{2\ell+1}$ & $\om_1$ & $2$ & $n_w$ ($\ell$ odd), $n_w\,\varpi_2^\vee(-1)$ ($\ell$ even)\\
$C_\ell$ & $Sp_{2\ell}$ & $\om_\ell$ & $2$ & $n_w$ ($\ell$ odd), $n_w\,\varpi_{\ell/2}^\vee(-1)$ ($\ell$ even)\\
$D_\ell$, $\ell$ odd & $Spin_{2\ell}$ & $\om_\ell$ & $4$ & $n_w$\\
$D_\ell$, $\ell$ even & $Spin_{2\ell}$ & $\om_1$, $\om_\ell$ & $2$, $2$ & $n_{w}$ for each\\
$E_6$ & $E_6^{sc}$ & $\om_1$ & $3$ & $n_w$\\
$E_7$ & $E_7^{sc}$ & $\om_7$ & $2$ & $n_w$
\end{tabular}
\medskip
\caption{Normalized lifts of generators of $\Om^{ad}$.}\label{tab:lifts}
\end{table}

\begin{proof}
For the classical types we use the standard coordinates $e_1,\dots,e_\ell$ of $\X_*(T)\otimes\Q$, and the dual coordinates $L_1,\dots,L_\ell$ of $\X^*(T)\otimes\Q$. We verify the criterion of Lemma~\ref{lem:criterion} type by type.

\emph{Type $A_\ell$.} Let $n=\ell+1$. Then $w$ is an $n$-cycle, and any lift $\dot x$ of $w$ is a monomial matrix of determinant $1$. Its $n$-th power is the product of its non-zero entries times $\Id_n$. Since $\det\dot x=1$ and the sign of an $n$-cycle is $(-1)^{n-1}$, we get $\dot x^{\,n}=(-1)^{n-1}\Id_n$. On the other hand $\exp(\varpi_1^\vee)=e^{-2\pi i/n}\Id_n$, so $\exp(\varpi_1^\vee)^{n(n-1)/2}=(-1)^{n-1}\Id_n$.

\emph{Type $B_\ell$.} Here $Q^\vee=\{\sum a_ie_i:\sum a_i\text{ even}\}$, $\varpi_1^\vee=e_1$ and $\varpi_2^\vee=e_1+e_2$. The element $w$ sends $e_1$ to $-e_1$ and fixes $e_2,\dots,e_\ell$. Since $w^2=1$, the set $I_1(w)$ consists of the positive roots made negative by $w$, namely $L_1$ and $L_1\pm L_j$ for $j>1$. So
\[
\gamma_w=2e_1+\sum_{j>1}\big((e_1-e_j)+(e_1+e_j)\big)=2\ell e_1.
\]
For $\lambda=0$, the criterion asks that $2(\ell-1)e_1\in2Q^\vee$, which holds if and only if $\ell$ is odd. For $\lambda=e_1+e_2$ we have $\lambda+w\lambda=2e_2$. The criterion asks that $2\big((\ell-1)e_1+e_2\big)\in2Q^\vee$, which holds if and only if $\ell$ is even.

\emph{Type $C_\ell$.} Here $Q^\vee=\bigoplus_i\Z e_i$, $\varpi_\ell^\vee=\frac12\sum_ie_i$ and $\varpi_{\ell/2}^\vee=e_1+\dots+e_{\ell/2}$. The element $w$ sends $e_i$ to $-e_{\ell+1-i}$. The positive roots made negative by $w$ are $L_i+L_j$ ($i<j$) and $2L_i$. So
\[
\gamma_w=\sum_{i<j}(e_i+e_j)+\sum_ie_i=\ell\sum_ie_i.
\]
For $\lambda=0$, the criterion asks that $(\ell-1)\sum_ie_i\in2Q^\vee$, which holds if and only if $\ell$ is odd. For $\lambda=e_1+\dots+e_{\ell/2}$, the vector $\gamma_w+\lambda+w\lambda-2\varpi_\ell^\vee$ has coefficient $\ell$ at $e_i$ for $i\leq\ell/2$, and $\ell-2$ for $i>\ell/2$. Both are even when $\ell$ is even.

\emph{Type $D_\ell$, $\ell$ odd.} Here $Q^\vee=\{\sum a_ie_i:\sum a_i\text{ even}\}$ and $\varpi_\ell^\vee=\frac12\sum_ie_i$. The element $w$ sends $e_1\mapsto e_\ell$ and $e_i\mapsto-e_{\ell+1-i}$ for $i\geq2$, and it has order $4$. A direct computation with Lemma~\ref{lem:adrian} gives
\[
\gamma_w=2(\ell-1)e_1+2(\ell-2)(e_2+\dots+e_{\ell-1}).
\]
Hence
\[
\gamma_w-12\varpi_\ell^\vee=2\big((\ell-4)e_1+(\ell-5)(e_2+\dots+e_{\ell-1})-3e_\ell\big).
\]
The sum of the coefficients inside the brackets is $(\ell-4)+(\ell-2)(\ell-5)-3$, which is even. So the criterion holds with $\lambda=0$.

\emph{Type $D_\ell$, $\ell$ even.} Let $w_1$ and $w_\ell$ be the elements for $\om_1$ and $\om_\ell$. Then $w_1$ sends $e_1\mapsto-e_1$, $e_\ell\mapsto-e_\ell$ and fixes the other $e_i$, while $w_\ell$ sends $e_i\mapsto-e_{\ell+1-i}$. A direct computation gives $\gamma_{w_1}=2(\ell-1)e_1$ and $\gamma_{w_\ell}=(\ell-1)\sum_ie_i$. So
\[
\gamma_{w_1}-2\varpi_1^\vee=2(\ell-2)e_1,\qquad\gamma_{w_\ell}-2\varpi_\ell^\vee=2\cdot\tfrac{\ell-2}{2}\sum_ie_i.
\]
Both lie in $2Q^\vee$ because $\ell$ is even. For the commutation, Lemma~\ref{lem:adrian}(1) gives
\[
\sum_{\alpha\in F(w_1,w_\ell)}\alpha^\vee=\sum_ie_i+(\ell-2)e_\ell,\qquad
\sum_{\alpha\in F(w_\ell,w_1)}\alpha^\vee=(\ell-1)e_1-(e_2+\dots+e_{\ell-1})+e_\ell.
\]
Their difference is $2\big(\tfrac{2-\ell}{2}e_1+(e_2+\dots+e_{\ell-1})+\tfrac{\ell-2}{2}e_\ell\big)$, which lies in $2Q^\vee$. Since $w_1w_\ell=w_\ell w_1$, we conclude that $n_{w_1}n_{w_\ell}=n_{w_\ell}n_{w_1}$.

\emph{Type $E_6$.} A direct computation gives
\[
\gamma_w=16\alpha_1^\vee+16\alpha_2^\vee+24\alpha_3^\vee+32\alpha_4^\vee+24\alpha_5^\vee+16\alpha_6^\vee\in2Q^\vee.
\]
Since $3\varpi_1^\vee\in Q^\vee$, also $6\varpi_1^\vee\in2Q^\vee$. So the criterion holds with $\lambda=0$.

\emph{Type $E_7$.} A direct computation gives
\[
\gamma_w=18\alpha_1^\vee+27\alpha_2^\vee+36\alpha_3^\vee+54\alpha_4^\vee+45\alpha_5^\vee+36\alpha_6^\vee+27\alpha_7^\vee=18\varpi_7^\vee.
\]
So $\gamma_w-2\varpi_7^\vee=8\cdot2\varpi_7^\vee$, which lies in $2Q^\vee$ since $2\varpi_7^\vee\in Q^\vee$.
\end{proof}

\subsection{Normalized lifts of all elements}
We also need normalized lifts of all elements of $\Om^{ad}$, not only of generators. Powers of a normalized lift are not always normalized (Example~\ref{exmp:all-elements}). The following lemma corrects this by a central element.

\begin{lem}\label{lem:powers}
Suppose $\Om^{ad}$ is cyclic of order $r$, generated by $\om_0=w\,t^{\mu_0}$ with a normalized lift $\dot x$. If $r$ is even, let $c=\exp(\mu_0)^{r/2}$. For $0\leq a<r$, let $d=\gcd(a,r)$ and put
\[
\dot x_a=\begin{cases}\dot x^{\,a}&\text{if }d\text{ is odd},\\ \dot x^{\,a}\,c^{a/d}&\text{if }d\text{ is even}.\end{cases}
\]
Then $\dot x_a$ is a normalized lift of $\om_0^a$.
\end{lem}

\begin{proof}
Let $s=r/d$ be the order of $\om_0^a$, and write $a=da'$ and $\om_0^a=w^at^{\mu}$. Then $\exp(\mu)=\exp(\mu_0)^a$, since $\mu$ is a sum of $a$ elements of $W\mu_0$ and $W$ acts trivially on $P^\vee/Q^\vee$. Put $z_0=\exp(\mu_0)$, a central element of order $r$. Let $b=0$ if $d$ is odd and $b=a'r/2$ if $d$ is even, so that $\dot x_a=\dot x^{\,a}z_0^{\,b}$. Since $as=a'r$, we get
\[
\dot x_a^{\,s}=(\dot x^{\,r})^{a'}\,z_0^{\,bs}=z_0^{\,M},\qquad M=a'\frac{r(r-1)}2+bs.
\]
We must show that $M\equiv M'$ modulo $r$, where $M'=a\frac{s(s-1)}2$. A direct computation gives
\[
M'-a'\frac{r(r-1)}2=a'ds^2\cdot\frac{1-d}2.
\]
If $d$ is odd, then $b=0$, and this is a multiple of $ds=r$. If $d$ is even, then $bs=a'ds^2/2$, and
\[
M'-M=-a'd^2s^2/2=-r\cdot\frac{a'ds}2,
\]
which is a multiple of $r$ because $d$ is even. Finally, Definition~\ref{defn:normalized} requires $\dot x_a\in N_G(T)(\Z)$. This holds because $z_0^{\,b}$ is $1$ or $c^{a'}$, and $c\in Z(G)(\Z)$ by \eqref{eq:exp-power}.
\end{proof}

\begin{cor}\label{cor:family}
There is an abelian subgroup $E\subset N_G(T)(\Z)$ and, for every $\om\in\Om^{ad}$, a normalized lift $\dot x_\om\in E$ of $\om$. If $G$ is not of type $B_{2m}$ or $C_{2m}$, then every $\dot x_\om$ is a product of powers of the Tits lifts in Table~\ref{tab:lifts} and an element of $Z(G)(\Z)$.
\end{cor}

\begin{proof}
If $\Om^{ad}$ is cyclic, let $\dot x$ be the lift in Table~\ref{tab:lifts}, and let $E$ be generated by $\dot x$ and $Z(G)(\Z)$. Lemma~\ref{lem:powers} gives the lifts, since $c\in Z(G)(\Z)$ by \eqref{eq:exp-power}.

In type $D_\ell$ with $\ell$ even, let $E$ be generated by $n_{w_1}$, $n_{w_\ell}$ and $Z(G)(\Z)$. It is abelian by Proposition~\ref{prop:generators}. We take the lifts $1$, $n_{w_1}$, $n_{w_\ell}$ and $n_{w_1}n_{w_\ell}$. The last one is normalized, because
\[
(n_{w_1}n_{w_\ell})^2=n_{w_1}^2n_{w_\ell}^2=\exp(\varpi_1^\vee+\varpi_\ell^\vee).\qedhere
\]
\end{proof}

\section{Splitting the Kottwitz homomorphism}\label{sec:thmA}

We now glue the normalized lifts of Section~\ref{sec:normalized} along the center, and prove Theorem~\ref{thm:A}.

The key point is that the normalization \eqref{eq:normalization} is a condition on each simple factor separately, and that it does not mention the subgroup $K$ of the center by which the factors are glued. If we raise the normalized lifts in all factors to a common power $r$, then all the $r$-th powers are given by the same formula $\exp(\mu_i)^{r(r-1)/2}$. So their product is a power of $\exp(\mu)$, where $\mu$ is the translation part, and $\exp(\mu)$ lies in $K$. This handles semisimple groups (Proposition~\ref{prop:semisimple}). In terms of Lemma~\ref{lem:cocycle}, it says that cocycles are stable under products and inflation (Remark~\ref{rem:cocycle-gluing}). The reductive case reduces to the semisimple one, because the torsion part of $\Om$ comes from the derived group.

\subsection{Assembling lifts}
We first record how commuting lifts of generators give a homomorphism.

\begin{lem}\label{lem:assemble}
Let $\Om=\bigoplus_k\langle\eta_k\rangle$ be a decomposition into cyclic subgroups, and write $\eta_k=w_k\,t^{\mu_k}$. Let $y_k\in N_G(T)(\Z)$ be lifts of $w_k$ that commute with each other. Suppose that $y_k^{r_k}=1$ whenever $\eta_k$ has finite order $r_k$. Then there is a unique homomorphism
\[
\psi\colon\Om\ra N_G(T)(\Z[t,t^{-1}])\quad\text{with}\quad\psi(\eta_k)=y_k\,t^{\mu_k},
\]
and $\psi(\om)$ maps to $\om$ in $\extw$ for all $\om$.
\end{lem}

\begin{proof}
We show that the elements $y_kt^{\mu_k}$ commute and have the right orders. For two of them we have
\[
(y_jt^{\mu_j})(y_kt^{\mu_k})=y_jy_k\,t^{w_k^{-1}\mu_j+\mu_k},\qquad (y_kt^{\mu_k})(y_jt^{\mu_j})=y_ky_j\,t^{w_j^{-1}\mu_k+\mu_j}.
\]
The $y$'s commute, and the exponents agree because $\eta_j\eta_k=\eta_k\eta_j$ in $\extw$. Similarly, $(y_kt^{\mu_k})^{r}=y_k^{\,r}\,t^{\nu}$, where $t^{\nu}$ is the translation part of $\eta_k^{\,r}$. If $\eta_k$ has order $r_k$, then $\nu=0$ for $r=r_k$, and $(y_kt^{\mu_k})^{r_k}=y_k^{r_k}=1$. The last claim follows from Lemma~\ref{lem:integral-ext}.
\end{proof}

\subsection{Semisimple groups}
Here is where the cocycle condition is used.

\begin{prop}\label{prop:semisimple}
Suppose $G$ is semisimple. For every $\eta\in\Om$ there is a lift $y_\eta\in N_G(T)(\Z)$ of $\pr(\eta)$ with $y_\eta^{\,r}=1$, where $r$ is the order of $\eta$. All the $y_\eta$ lie in one abelian subgroup of $N_G(T)(\Z)$.
\end{prop}

\begin{proof}
Write $G=G^{sc}/K$ with $G^{sc}=\prod_iG_i$, where each $G_i$ is almost simple and simply connected. Then $\Om\subset\Om^{ad}=\prod_i\Om_{G_i^{ad}}$. For each $i$, choose an abelian group $E_i$ and normalized lifts $\dot x_\om\in E_i$ for $\om\in\Om_{G_i^{ad}}$, as in Corollary~\ref{cor:family}.

Let $\eta\in\Om$ have order $r$, and write its components as $\eta_i=w_i\,t^{\mu_i}\in\Om_{G_i^{ad}}$. Put
\[
\dot x=\prod_i\dot x_{\eta_i}\in G^{sc}(\Z),\qquad y_\eta=\pi(\dot x)\in N_G(T)(\Z).
\]
Then $y_\eta$ lifts $\pr(\eta)$, and it lies in the abelian group $\pi(\prod_iE_i)$. It remains to show that $y_\eta^{\,r}=1$, that is, $\dot x^{\,r}\in K$. We check this in $G^{sc}(\CC)$.

Let $r_i$ be the order of $\eta_i$, and write $r=m_ir_i$. By Lemma~\ref{lem:exp}, the element $\exp(\mu_i)$ has order $r_i$. Since $\dot x_{\eta_i}$ is a normalized lift of $\eta_i$,
\[
\dot x_{\eta_i}^{\,r}=\big(\dot x_{\eta_i}^{\,r_i}\big)^{m_i}=\exp(\mu_i)^{m_ir_i(r_i-1)/2}=\exp(\mu_i)^{r(r_i-1)/2}.
\]
We claim that this equals $\exp(\mu_i)^{r(r-1)/2}$. Indeed, the two exponents differ by
\[
\frac{r(r-1)}2-\frac{r(r_i-1)}2=\frac{r(r-r_i)}2=r_i^2\cdot\frac{m_i(m_i-1)}2,
\]
which is a multiple of $r_i$. The elements $\dot x_{\eta_i}$ lie in different simple factors, so they commute. Hence
\[
\dot x^{\,r}=\prod_i\exp(\mu_i)^{r(r-1)/2}=\exp(\mu)^{r(r-1)/2},\qquad\text{where }\mu=\sum_i\mu_i.
\]
Here $\mu$ is the translation part of $\eta$. It lies in $\X_*(T)$ because $\eta\in\Om$, so $\exp(\mu)\in K(\CC)$ by Lemma~\ref{lem:exp}(2). Therefore $\dot x^{\,r}\in K(\CC)$, and $y_\eta^{\,r}=1$.
\end{proof}

The point of the proof is that, after raising to the common power $r$, the normalization gives the same exponent $r(r-1)/2$ in every factor. Here is the conceptual reason.

\begin{rem}[the cocycle point of view]\label{rem:cocycle-gluing}
We keep the notation of the proof and work over $\CC$. Put
\[
g_i=\dot x_{\eta_i}\,t^{\mu_i}.
\]
Since $r_i\mu_i\in Q^\vee$, this element lies in $G_i(\cK_{r_i})$. By Lemma~\ref{lem:cocycle}, the normalization of $\dot x_{\eta_i}$ says exactly that $g_i$ is a $1$-cocycle of $\Gamma_{r_i}$, that is, $N_{r_i}(g_i)=1$.

Now regard $g_i$ as an element of $G_i(\cK_r)$. By Section~\ref{subsec:galois}, the generator $\gamma$ of $\Gamma_r$ acts on $\cK_{r_i}$ as the generator of $\Gamma_{r_i}$, and $\gamma^{r_i}$ acts trivially on $\cK_{r_i}$. So the sequence $g_i,\gamma(g_i),\gamma^2(g_i),\dots$ is periodic with period $r_i$. Grouping the $r$ factors of $N_r(g_i)$ into $m_i$ consecutive blocks of length $r_i$, we get
\[
N_r(g_i)=N_{r_i}(g_i)^{m_i}=1.
\]
This is the inflation of a cocycle from $\Gamma_{r_i}$ to $\Gamma_r$. Applying Lemma~\ref{lem:cocycle} with $r$ in place of $r_i$, we see that this is equivalent to $\dot x_{\eta_i}^{\,r}=\exp(\mu_i)^{r(r-1)/2}$. This is the claim in the proof.

Finally, $g=\prod_ig_i=\dot x\,t^{\mu}$ satisfies $N_r(g)=\prod_iN_r(g_i)=1$, because the $g_i$ commute and $\gamma$ preserves each factor. So $g$ is a cocycle, and Lemma~\ref{lem:cocycle} gives $\dot x^{\,r}=\exp(\mu)^{-r(r-1)/2}$, which lies in $K$. In short: cocycle lifts are stable under inflation and products, and the $r$-th power of a cocycle lift is a power of its monodromy $\exp(\mu)$, which lies in $K$.
\end{rem}

\subsection{Proof of Theorem~\ref{thm:A}}
The reductive case follows by splitting $\Om$ into its torsion part, which comes from the derived group, and a free part.

\begin{thm}\label{thm:A-body}
Let $G$ be a connected reductive group split over $\Z$. There is a homomorphism $\psi\colon\Om\ra N_G(T)(\Z[t,t^{-1}])$ such that $\psi(\om)$ maps to $\om$ in $\extw$. For every field $k$, the induced map $\Om\ra LG(k)$ is a homomorphic section of $\kappa_G$, and
\[
LG_k\cong LG_k^{\circ}\rtimes\Om.
\]
\end{thm}

\begin{proof}
Let $M=G_{\der}$ and $T_M=T\cap M$. The lattice $\X_*(T_M)$ is saturated in $\X_*(T)$, because $T_M$ is a subtorus. So the torsion subgroup of $\Om=\X_*(T)/Q^\vee$ is $\X_*(T_M)/Q^\vee=\Om_M$. Choose a decomposition $\Om=\Om_M\oplus F$ with $F$ free, and a basis $f_1,\dots,f_m$ of $F$. Decompose $\Om_M$ into cyclic subgroups $\langle\eta_k\rangle$. For each $\eta_k$, let $y_k$ be the lift from Proposition~\ref{prop:semisimple} applied to $M$.

It remains to lift the basis vectors $f_j$ of $F$. Let $\bar f_j$ be the image of $f_j$ in $\Om_{G^{ad}}=\Om_{M^{ad}}$. The map $\extw_G\ra\extw_{G^{ad}}$ is the identity on $W$, so $f_j$ and $\bar f_j$ have the same finite part $w_j$. Write $M^{sc}=\prod_iG_i$ and $\bar f_j=(\bar f_{j,i})_i$ with $\bar f_{j,i}\in\Om_{G_i^{ad}}$, and put
\[
y_j'=\pi_M\Big(\prod_i\dot x_{\bar f_{j,i}}\Big),
\]
where $\pi_M\colon M^{sc}\ra M$ is the simply connected cover and the $\dot x$ are the normalized lifts used in the proof of Proposition~\ref{prop:semisimple}. Then $y_j'$ lifts $w_j$. Since $f_j$ has infinite order, Lemma~\ref{lem:assemble} imposes no condition on the order of $y_j'$. All the elements $y_k$ and $y_j'$ lie in the abelian group $\pi_M(\prod_iE_i)\subset N_M(T_M)(\Z)\subset N_G(T)(\Z)$, so Lemma~\ref{lem:assemble} gives $\psi$.

For a field $k$, the image of $\tau n_wt^{\mu}$ in $LG(k)$ has Kottwitz invariant $[\mu]$. So $\kappa_G\circ\psi=\id_{\Om}$. Since $LG_k$ is the disjoint union of the cosets $LG_k^{\circ}\psi(\om)$, the map $(g,\om)\mapsto g\psi(\om)$ is an isomorphism $LG_k^{\circ}\rtimes\Om\cong LG_k$.
\end{proof}

Since $\psi$ is defined over $\Z[t,t^{-1}]$, we may also substitute a uniformizer for $t$. This gives the same result for $p$-adic groups, for which no loop group is available.

\begin{cor}\label{cor:local-field}
Let $F$ be a non-archimedean local field with ring of integers $\cO_F$ and uniformizer $\varpi_F$. Let $\kappa_G\colon G(F)\ra\Om$ be the Kottwitz homomorphism, normalized by $\kappa_G(\mu(\varpi_F))=[\mu]$. Then the composite of $\psi$ with the map $N_G(T)(\Z[t,t^{-1}])\ra G(F)$, $t\mapsto\varpi_F$, is a homomorphic section of $\kappa_G$.
\end{cor}

\begin{proof}
The image of $\tau n_w$ lies in $G(\cO_F)$, on which $\kappa_G$ is trivial \cite[\S7]{Ko97}, \cite{HR08}. The image of $t^{\mu}$ is $\mu(\varpi_F)$.
\end{proof}

\subsection{The normalization is necessary}\label{subsec:necessary}
We end this section with two examples. The first shows that the normalization is forced as soon as a factor of type $A$ is involved. The second shows that we need it for all elements of $\Om^{ad}$, and not only for generators.

\begin{exmp}\label{exmp:forced}
Let $H$ be almost simple and simply connected, and let $\om_H=w_Ht^{\mu_H}\in\Om_{H^{ad}}$ have order $r$. Let $\om_A=w_At^{\mu_A}$ generate $\Om_{PGL_r}$. Let $K\subset Z(SL_r\times H)$ be generated by $(\exp\mu_A,\exp\mu_H)$, and put $G=(SL_r\times H)/K$. Then $\Om\cong\Z/r\Z$ is generated by $\eta=(\om_A,\om_H)$.

A lift of $\pr(\eta)$ to $N_G(T)$ is the image of a pair $(x,y)$. By Proposition~\ref{prop:generators}, $x^r=\exp(\mu_A)^{N}$ with $N=r(r-1)/2$. So the lift has order $r$ if and only if $(x^r,y^r)\in K$, that is, if and only if $y^r=\exp(\mu_H)^N$. In other words, $y$ must be a normalized lift of $\om_H$. Taking $H=Sp_{2n}$, $E_6$, $Spin_6$ or $Spin_8$ gives examples with $r=2,3,4,2$. A similar example with $G^{sc}=SL_2^3\times Spin_8$ explains why we need commuting lifts in type $D_\ell$ with $\ell$ even.
\end{exmp}

\begin{exmp}\label{exmp:all-elements}
Let $G=(SL_6\times SL_3)/\mu_3$, with $\mu_3$ embedded diagonally. Then $\Om\cong\Z/3\Z$. A generator is $\eta=(\om^2,\om')$, where $\om$ generates $\Om_{PGL_6}\cong\Z/6\Z$ and $\om'$ generates $\Om_{PGL_3}\cong\Z/3\Z$. The finite parts of $\om$ and $\om'$ are a $6$-cycle and a $3$-cycle. Let $\dot x\in SL_6$ and $\dot y\in SL_3$ be lifts of them. They are normalized lifts, with
\[
\dot x^{\,6}=-\Id_6,\qquad\dot y^{\,3}=\Id_3.
\]
The obvious lift of $\eta$ is the image of $(\dot x^{\,2},\dot y)$. But $(\dot x^{\,2},\dot y)^3=(-\Id_6,\Id_3)$, which does not lie in $\mu_3$. So this lift has order $6$, not $3$.

The reason is that $\dot x^{\,2}$ is not a normalized lift of $\om^2$. Indeed, $\om^2$ has order $3$, so a normalized lift of $\om^2$ must have cube $\exp(2\varpi_1^\vee)^3=\Id_6$. Lemma~\ref{lem:powers} gives the normalized lift $-\dot x^{\,2}$, and the image of $(-\dot x^{\,2},\dot y)$ has order $3$.
\end{exmp}

\section{Pinned splittings}\label{sec:pinned}

A splitting of $\kappa_G$ is \emph{pinned} if its image preserves an affine pinning. In this section we show that Tits' lifts give pinned splittings in many cases (Theorem~\ref{thm:B}(1)), and that pinned splittings do not exist for two families of groups (Theorem~\ref{thm:B}(2)).

The difference between the two situations is the parity phenomenon in types $B_{2m}$ and $C_{2m}$ found in Proposition~\ref{prop:generators}. In those types, Tits' lift is not normalized. The correction that normalizes it is an element of $T(\Z)$, and it does not preserve the pinning. When such a factor is glued to a factor of type $A$ along the center, the normalization is forced (Example~\ref{exmp:forced}), and no pinned splitting exists (Proposition~\ref{prop:B2-body}).

\subsection{Tits' lifts preserve an affine pinning}
This subsection contains the main step in the proof of Theorem~\ref{thm:B}(1): Tits' lifts preserve a suitable affine pinning (Lemma~\ref{lem:pinned-tits}).

\begin{defn}\label{defn:P-omega}
Let $G$ be almost simple and simply connected, and let $\om=w\,t^{\varpi^\vee}\in\Om^{ad}$ be non-trivial. Put $g_\om=n_w\,\varpi^\vee(t)\in N_{G^{ad}}(T^{ad})(\Z[t,t^{-1}])$, $a_0=\delta-\theta$ and $\beta=\om(a_0)\in\Delta$. Let $\cP_\om$ be the affine pinning of $LG$ with root vectors
\[
x_\alpha=u_\alpha(1)\ \ (\alpha\in\Delta),\qquad x_{a_0}=g_\om^{-1}\,u_\beta(1)\,g_\om.
\]
\end{defn}

Here $n_w$ denotes the image of the Tits lift in $G^{ad}$, and $\beta$ is a simple root because only the trivial element of $\Om^{ad}$ fixes $a_0$. Since $n_w$ maps root vectors to root vectors up to sign, $x_{a_0}=u_{-\theta}(\pm t)$; so $\cP_\om$ is defined over $\Z$.

\begin{lem}\label{lem:pinned-tits}
In the setting of Definition~\ref{defn:P-omega}, suppose that $n_w^r$ is central, where $r$ is the order of $\om$. Then $g_\om$ preserves $\cP_\om$.
\end{lem}

\begin{proof}
Write $g=g_\om$. Let $\alpha\in\Delta$ with $\om(\alpha)\in\Delta$. Then $\langle\alpha,\varpi^\vee\rangle=0$ and $w(\alpha)=\om(\alpha)$. So $g\,x_\alpha\,g^{-1}=n_wu_\alpha(1)n_w^{-1}$, and this equals $x_{\om(\alpha)}$ by Lemma~\ref{lem:springer}.

Consider the orbit of $a_0$, that is, $a_0,\om(a_0),\dots,\om^{r-1}(a_0)$. Its elements other than $a_0$ lie in $\Delta$. By definition $g\,x_{a_0}g^{-1}=x_\beta$, and the steps inside $\Delta$ are handled by Lemma~\ref{lem:springer}. For the last step, we compute $g^r$ in $LG^{ad}$. Since $\om^r=1$, the translation parts cancel, and $g^r$ is the image of $n_w^r\in G$ under the projection $G\ra G^{ad}$. This image is trivial because $n_w^r$ is central. Hence $g^r=1$, and
\[
g\,x_{\om^{r-1}(a_0)}\,g^{-1}=g^{r-1}x_\beta\,g^{1-r}=g^r\,x_{a_0}\,g^{-r}=x_{a_0}.
\]
Every other orbit lies in $\Delta$.
\end{proof}

All lifts $n_w$ in Table~\ref{tab:lifts}, and also the uncorrected Tits lifts in types $B_{2m}$ and $C_{2m}$, have central $r$-th powers. In the latter case $n_w^2=1$ by the proof of Proposition~\ref{prop:generators}.

In type $D_\ell$ with $\ell$ even, the group $\Om^{ad}$ has two generators. The next lemma shows that one affine pinning works for both.

\begin{lem}\label{lem:pinned-D}
In type $D_\ell$ with $\ell$ even, let $g_1=g_{\om_1}$ and $g_\ell=g_{\om_\ell}$ be the elements of Definition~\ref{defn:P-omega}. Then $g_\ell$ also preserves $\cP_{\om_1}$.
\end{lem}

\begin{proof}
The elements $g_1$ and $g_\ell$ commute, by Proposition~\ref{prop:generators} and the computation in Lemma~\ref{lem:assemble}. Let $\beta_1=\om_1(a_0)$. Then
\[
g_\ell\,x_{a_0}\,g_\ell^{-1}=g_1^{-1}\big(g_\ell\,x_{\beta_1}g_\ell^{-1}\big)g_1=g_1^{-1}x_{\om_\ell(\beta_1)}g_1=x_{\om_1^{-1}\om_\ell\om_1(a_0)}=x_{\om_\ell(a_0)}.
\]
Here we used Lemma~\ref{lem:springer} twice, since $\om_\ell(\beta_1)$ and $\om_\ell(a_0)$ lie in $\Delta$. The roots in $\Delta$ are handled by Lemma~\ref{lem:springer}.
\end{proof}

As a first consequence, the extension \eqref{eq:pinned-ext} is exact for every semisimple group, in every characteristic.

\begin{cor}\label{cor:surjective}
Let $G$ be semisimple.
\begin{enumerate}[(1)]
\item There is an affine pinning $\cP$ of $LG$, defined over $\Z$, such that $\kappa_G(\tOm_\cP)=\Om$ over every field $k$.
\item If $k$ is algebraically closed, then $\kappa_G(\tOm_{\cP'})=\Om$ for every affine pinning $\cP'$ of $LG_k$.
\end{enumerate}
In particular, the sequence \eqref{eq:pinned-ext} is exact in both cases.
\end{cor}

\begin{proof}
(1) Write $G=G^{sc}/K$ with $G^{sc}=\prod_iG_i$. Let $\cP$ be the product of the affine pinnings $\cP_\om$ of Definition~\ref{defn:P-omega}, one for each factor, with $\om=\om_1$ in type $D_\ell$ for $\ell$ even. Let $\eta\in\Om$, with translation part $t^{\mu}$. Each component of $\eta$ is a product of powers of the generators in Table~\ref{tab:lifts}. Let $\dot x_i$ be the same product of powers of the corresponding Tits lifts, and put
\[
g_\eta=\pi\Big(\prod_i\dot x_i\Big)\,t^{\mu}\in N_G(T)(\Z[t,t^{-1}]).
\]
This element lies over $\eta$. Its image in $LG^{ad}=\prod_iLG_i^{ad}$ is a product of powers of the elements $g_\om$ of Definition~\ref{defn:P-omega}. So it preserves $\cP$, by Lemmas~\ref{lem:pinned-tits} and~\ref{lem:pinned-D}. Here we also use the uncorrected Tits lifts in types $B_{2m}$ and $C_{2m}$, whose squares are trivial.

(2) This is \cite[Proposition~3.14]{DLYZ} in the untwisted case. There, the characteristic of $k$ is assumed to be sufficiently large; in fact, the proof works in every characteristic. For completeness, we give the proof here.

By Lemma~\ref{lem:conj-pinning}, we may assume that $\cP'=\cP_c$ for some $c\in(k^\times)^{\Daff}$, with $\cP$ as in (1). Fix $\eta\in\Om$. The element $g_\eta$ of (1) preserves $\cP$, so it maps $u_a(x)$ to $u_{\eta(a)}(x)$ for all $a\in\Daff$ and $x\in k$. Hence, for $\tau\in T(k)$, the element $g_\eta\tau$ maps $u_a(c_a)$ to $u_{\eta(a)}(a(\tau)\,c_a)$. So $g_\eta\tau$ preserves $\cP_c$ if and only if
\[
a(\tau)=c_{\eta(a)}/c_a\qquad\text{for all }a\in\Daff.\tag{$*$}
\]
We show that $(*)$ has a solution. The characters $\pi(a)$ of $T^{ad}$, for $a\in\Daff$, generate $\X^*(T^{ad})$. Their only relations are $\sum_an_a\pi(a)=0$, one for each simple factor, where $a$ runs over the simple affine roots of that factor and the $n_a$ are the labels of the affine Dynkin diagram. So $(*)$ has a solution $\bar\tau\in T^{ad}(k)$ if and only if $\prod_a(c_{\eta(a)}/c_a)^{n_a}=1$ for each simple factor. This holds because $\eta$ permutes the simple affine roots of each factor and preserves their labels. Since $k$ is algebraically closed, $\bar\tau$ lifts to $\tau\in T(k)$, and $g_\eta\tau\in\tOm_{\cP_c}$ lies over $\eta$.
\end{proof}

\subsection{Proof of Theorem~\ref{thm:B}(1)}
We combine the two lemmas with the construction of Section~\ref{sec:thmA}.

\begin{thm}\label{thm:B1-body}
Suppose that $G$ is almost simple, or that no simple factor of $G$ is of type $B_{2m}$ or $C_{2m}$. Then $G$ is pinnable over every field $k$. If no simple factor is of type $B_{2m}$ or $C_{2m}$, the section $\psi$ of Theorem~\ref{thm:A-body} can be chosen to preserve an affine pinning defined over $\Z$.
\end{thm}

\begin{proof}
First suppose that no simple factor of $G$ is of type $B_{2m}$ or $C_{2m}$. Let $M^{sc}=\prod_iG_i$ be the simply connected cover of $G_{\der}$. For each $i$, let $\cP_i$ be the affine pinning of Definition~\ref{defn:P-omega} attached to the generator of $\Om_{G_i^{ad}}$ in Table~\ref{tab:lifts}. In type $D_\ell$ with $\ell$ even, we use $\om_1$. Together they form an affine pinning $\cP$ of $LG$, defined over $\Z$.

The action of $LG$ on affine pinnings factors through $LG^{ad}=\prod_iLG_i^{ad}$. In the proof of Theorem~\ref{thm:A-body}, every $\psi(\om)$ is built from the lifts of Corollary~\ref{cor:family}. These are products of powers of Tits lifts and central elements. Hence the image of $\psi(\om)$ in $LG_i^{ad}$ is a product of powers of the elements $g_\om$ of Definition~\ref{defn:P-omega}. By Lemmas~\ref{lem:pinned-tits} and~\ref{lem:pinned-D}, $\psi(\om)$ preserves $\cP$.

Now suppose $G$ is almost simple of type $B_{2m}$ or $C_{2m}$. Then $G$ is simply connected or adjoint. If $G$ is simply connected, then $\Om$ is trivial. If $G$ is adjoint, the element $g_\om$ of Definition~\ref{defn:P-omega} lies in $N_G(T)(\Z[t,t^{-1}])$. It preserves $\cP_\om$ by Lemma~\ref{lem:pinned-tits}, and $g_\om^2=1$ because $n_w^2=1$. So $\om\mapsto g_\om$ is a pinned splitting.
\end{proof}

\subsection{Non-pinnable groups}
We now show that the two families of Theorem~\ref{thm:B}(2) are not pinnable. We prove more: for every affine pinning, every lift of the generator in the stabilizer has order $4$. The argument compares an arbitrary affine pinning with a standard one, using Lemma~\ref{lem:conj-pinning}.

\begin{prop}\label{prop:B2-body}
Let $k$ be a field with $\chara k\neq2$, and let $n\geq2$ be even. Let $H$ be $Sp_{2n}$ or $Spin_{2n+1}$, let $z$ be the non-trivial element of $Z(H)(k)$, and put
\[
G=(SL_2\times H)/\langle(-\Id,z)\rangle.
\]
Then $\Om\cong Z(G)(k)\cong\Z/2\Z$. For every affine pinning $\cP$, every element of $\tOm_{\cP}$ that maps to the generator of $\Om$ has order $4$. For the affine pinning built from Definition~\ref{defn:P-omega} in both factors, $\tOm_{\cP}\cong\Z/4\Z$. In particular, $G$ is not pinnable over $k$.
\end{prop}

\begin{proof}
The group $\Om$ is generated by $\eta=(\om_A,\om_H)$, where $\om_A$ and $\om_H$ are the non-trivial elements of $\Om_{PGL_2}$ and $\Om_{H^{ad}}$. Its finite part $w=\pr(\eta)$ is an involution.

\emph{The standard pinning.} Let $g$ be the image of $(\sigma t^{\mu_A},n_wt^{\mu_H})$, where $\sigma$ and $n_w$ are the Tits lifts in $SL_2$ and $H$. Then $g\in LG(k)$, and $g$ preserves the affine pinning $\cP$ built from Definition~\ref{defn:P-omega}, by Lemma~\ref{lem:pinned-tits}. We have $\sigma^2=-\Id$, and $n_w^2=1$ because $n$ is even (see the proof of Proposition~\ref{prop:generators}). So $g^2=z_0$, the image of $(-\Id,1)$. Since $\chara k\neq2$, this is the non-trivial element of $Z(G)(k)$. By Proposition~\ref{prop:Z-extension},
\[
\tOm_{\cP}=\{1,z_0,g,gz_0\}\cong\Z/4\Z.
\]

\emph{Every affine pinning.} Let $\cP'$ be any affine pinning, and let $g'\in\tOm_{\cP'}$ map to $\eta$. We show that $g'^2=z_0$, so that $g'$ has order $4$. We may assume that $k$ is algebraically closed. By Lemma~\ref{lem:conj-pinning}, applied to $g\in\tOm_\cP$ and $g'\in\tOm_{\cP'}$, we may replace $g'$ by a conjugate and assume that $g'=g\tau$ with $\tau\in T(k)$. Moreover, there is $c\in(k^\times)^{\Daff}$ with $a(\tau)=c_{\eta(a)}/c_a$ for all $a\in\Daff$. Since $\eta^2=1$, this gives
\[
a(\tau)=1\ \text{ if }\eta(a)=a,\qquad a(\tau)\,\eta(a)(\tau)=1\ \text{ for every }a\in\Daff.\tag{$*$}
\]

Next we compute $g'^2$. Since $g^2=z_0$ and $g^{-1}\tau g=w(\tau)$, we have $g'^2=g\tau g\tau=z_0\,w(\tau)\,\tau$. So for every $\lambda\in\X^*(T)$,
\[
\lambda(g'^2)=\lambda(z_0)\cdot(\lambda+w\lambda)(\tau).
\]
Since $g'^2\in Z(G)(k)=\{1,z_0\}$, it suffices to find $\lambda\in\X^*(T)$ with $\lambda(z_0)=-1$ and $(\lambda+w\lambda)(\tau)=1$. By $(*)$, the second condition holds if $\lambda+w\lambda$ is an integral combination of simple roots fixed by $\eta$ and of sums $\alpha+\eta(\alpha)$ of simple roots. We take $\lambda=\varpi_A+\lambda_H$, where $\varpi_A$ is the fundamental weight of $SL_2$ and $\lambda_H$ is as follows.
\begin{itemize}
\item If $H=Sp_{2n}$, take $\lambda_H=L_1$. Then
\[
\lambda+w\lambda=L_1-L_n=\alpha_{n/2}+\sum_{j<n/2}(\alpha_j+\alpha_{n-j}),
\]
and $\eta$ fixes $\alpha_{n/2}$ and swaps $\alpha_j$ and $\alpha_{n-j}$.
\item If $H=Spin_{2n+1}$, take $\lambda_H=\varpi_n=\frac12(L_1+\dots+L_n)$, the spin weight. Then
\[
\lambda+w\lambda=L_2+\dots+L_n=\sum_{j=2}^n(j-1)\,\alpha_j,
\]
and $\eta$ fixes $\alpha_2,\dots,\alpha_n$.
\end{itemize}
Here we used $\varpi_A+\sigma(\varpi_A)=0$. In both cases $\lambda_H(z)=-1$. So $\lambda$ is trivial on $(-\Id,z)$ and is a character of $T$, and $\lambda(z_0)=\varpi_A(-\Id)=-1$. Hence $\lambda(g'^2)=-1$, so $g'^2=z_0$ and $g'$ has order $4$. Therefore $G$ is not pinnable over $\bar k$, and hence not over $k$.
\end{proof}

When $k$ is algebraically closed, the case $H=Spin_{2n+1}$ also follows from the case $H=Sp_{2n}$ by Theorem~\ref{thm:dual}, since the two groups are Langlands dual to each other.

\begin{rem}[small characteristic]\label{rem:char2}
The assumption $\chara k\neq2$ in Proposition~\ref{prop:B2-body} is necessary. In characteristic $2$, the group $\mu_2(k)$ is trivial. In the proof of Proposition~\ref{prop:B2-body} we then have $g^2=1$, because $-\Id=\Id$. So $\tOm_{\cP}\cong\Om$, and $G$ is pinnable. At the same time, the kernel of $SL_2\ra PGL_2$ is infinitesimal, and the loop group of $PGL_2$ is not reduced \cite[Remarks~6.3 and~6.4]{PR08}.
\end{rem}

\section{Duality of the pinned extensions}\label{sec:duality}

In this section $G$ is a split semisimple group, and $G^\vee$ is its Langlands dual group, with dual torus $T^\vee$. We state and prove Theorem~\ref{thm:duality}, a more precise form of Theorem~\ref{thm:dual}.

We first prove the theorem over $\CC$. The idea is as follows. The stabilizer of an affine pinning fixes the sum of the root vectors of the pinning, an element of the loop Lie algebra. After conjugation by a suitable element of $T^{ad}$ over a ramified extension of $\cK$, this sum becomes Kostant's cyclic element, up to a scalar in each simple factor. Over $\CC$, Kostant's element is regular semisimple. So the conjugated stabilizer lies in a torus, and it is the group of fixed points of a Coxeter element on that torus (Proposition~\ref{prop:coxeter-torus}). These fixed points have a simple description in terms of lattices (Proposition~\ref{prop:lattice}), and the lattices for $G$ and for $G^\vee$ are dual to each other. The same arguments work over any algebraically closed field of characteristic $0$. In positive characteristic, Kostant's element need not be regular semisimple, and we lift pinned elements to characteristic $0$ instead (Section~\ref{subsec:char-p}). As a by-product, we get a lattice criterion for pinnability (Corollary~\ref{cor:criterion}).

\subsection{Statement}\label{subsec:dual-statement}
Let $k$ be an algebraically closed field, and let $\cP$ and $\cP^\vee$ be affine pinnings of $LG_k$ and $LG^\vee_k$. As in Section~\ref{subsec:affine-pinning}, we write $\tOm_\cP\subset LG(k)$ and $\tOm_{\cP^\vee}\subset LG^\vee(k)$ for their stabilizers.\footnote{In \cite{DLYZ}, these stabilizers are denoted $M_\phi$ and $M^\vee_{\phi^\vee}$.} By Proposition~\ref{prop:Z-extension} and Corollary~\ref{cor:surjective}, we have exact sequences
\[
1\ra Z(G)(k)\ra \tOm_\cP\ra\Om_G\ra1,\qquad 1\ra Z(G^\vee)(k)\ra \tOm_{\cP^\vee}\ra\Om_{G^\vee}\ra1,
\]
where $\Om_G=\X_*(T)/Q^\vee$ and $\Om_{G^\vee}=\X^*(T)/Q$. Here $Q$ is the root lattice of $G$, which is the coroot lattice of $G^\vee$. There are standard pairings
\[
Z(G)(k)\times\Om_{G^\vee}\ra k^\times,\ \ (z,[\lambda])\mapsto\lambda(z),
\]
\[
\Om_G\times Z(G^\vee)(k)\ra k^\times,\ \ ([\mu],z^\vee)\mapsto\mu(z^\vee),
\]
which are perfect when $\chara k\nmid|Z(G)|\cdot|\pi_1(G)|$.

\begin{thm}\label{thm:duality}
Let $k$ be an algebraically closed field, and let $c\in W$ be a Coxeter element.
\begin{enumerate}[(1)]
\item Suppose that $\chara k\nmid|Z(G)|$. Then the extension $1\ra Z(G)(k)\ra \tOm_\cP\ra\Om_G\ra1$ is isomorphic to
\[
0\ra Q^\vee/(1-c)\X_*(T)\ra\X_*(T)/(1-c)\X_*(T)\ra\X_*(T)/Q^\vee\ra0.
\]
In particular, $\tOm_\cP$ is commutative.
\item Suppose that $\chara k\nmid|Z(G)|\cdot|\pi_1(G)|$. Then there is a perfect pairing $\tOm_\cP\times \tOm_{\cP^\vee}\ra k^\times$ under which $Z(G)(k)$ and $Z(G^\vee)(k)$ are the annihilators of each other, and which induces the two standard pairings. In particular, the two extensions are Pontryagin dual to each other, and $G$ is pinnable over $k$ if and only if $G^\vee$ is.
\end{enumerate}
\end{thm}

We prove Theorem~\ref{thm:duality} for $k=\CC$ in Sections~\ref{subsec:kostant}--\ref{subsec:pairing}, and in general in Section~\ref{subsec:char-p}.

\subsection{Kostant's cyclic element}\label{subsec:kostant}
The main input is Kostant's theory of the cyclic element and the principal element \cite{Kos59}.

Write $G^{sc}=\prod_iG_i$ with $G_i$ almost simple. For each $i$, let $\Delta_i$, $\theta_i$, $h_i$ and $\rho_i^\vee$ be the simple roots, the highest root, the Coxeter number, and the half sum of the positive coroots of $G_i$. We write
\[
\rho^\vee/h:=\sum_i\rho_i^\vee/h_i\ \in\ \X_*(T^{ad})\otimes\Q.
\]
This is only a notation: there is no number $h$ unless $G$ is simple, in which case $\rho^\vee/h$ is $\rho^\vee$ divided by the Coxeter number $h$ of $G$. The \emph{principal element} is
\[
x_0=\exp(\rho^\vee/h)\in T^{ad}(\CC).
\]

Let $X_\beta\in\fg$ be root vectors. For constants $c_\beta\in\CC^\times$ put
\[
E_{0,i}=\sum_{\alpha\in\Delta_i}c_\alpha X_\alpha+c_{-\theta_i}X_{-\theta_i},\qquad E_0=\sum_iE_{0,i}.
\]
So $E_0$ is a sum of root vectors, one for each simple root and one for the negative of each highest root. It is called Kostant's \emph{cyclic element}.

\begin{fact}[{\cite[\S6]{Kos59}}]\label{fact:kostant}
\begin{enumerate}[(1)]
\item $E_0$ is regular semisimple. So $T_0:=Z_G(E_0)$ is a maximal torus of $G$.
\item $\Ad(x_0)E_{0,i}=e^{2\pi i/h_i}E_{0,i}$ for each $i$. Hence $x_0$ normalizes $T_0$, and it acts on $T_0$ by a Coxeter element of $W(G,T_0)$.
\end{enumerate}
We write $c$ for the action of $\Ad(x_0^{-1})$ on $T_0$. It is also a Coxeter element.
\end{fact}

Kostant states (1) for $c_\beta\equiv1$ and $G$ simple. The general case follows, because $E_{0,i}$ is conjugate under $T^{ad}_i(\CC)$ to a non-zero multiple of the standard element. The first claim in (2) is immediate: $\langle\alpha,\rho_i^\vee\rangle=1$ for $\alpha\in\Delta_i$, and $\langle-\theta_i,\rho_i^\vee\rangle=1-h_i$.

We also recall a well-known fact about Coxeter elements.

\begin{fact}[{\cite[Ch.~VI, \S1, Exercise~22]{Bou02}}]\label{fact:det}
For a Coxeter element $c$, we have
\[
|\det(1-c)|=|P^\vee/Q^\vee|\qquad\text{and}\qquad(c-1)P^\vee=Q^\vee.
\]
\end{fact}

\subsection{Pinned elements lie in a Coxeter torus}\label{subsec:coxeter-torus}
We now show that $\tOm_\cP$ is isomorphic to the group $T_0^{\,c}$ of fixed points of $c$ on $T_0$ (Proposition~\ref{prop:coxeter-torus}). The idea is that $\tOm_\cP$ fixes the sum of the root vectors of the affine pinning $\cP$, and that conjugation by the element $t^{\rho^\vee/h}$ defined below turns this sum into Kostant's cyclic element $E_0$, up to a scalar in each simple factor.

The root vectors of $\cP$ are $c_\alpha X_\alpha$ for $\alpha\in\Delta_i$ and $c_{-\theta_i}\,t\,X_{-\theta_i}$, for some constants $c_\beta\in\CC^\times$. For each $i$, let $E_i$ be the sum of the root vectors of $\cP$ that belong to the $i$-th factor:
\[
E_i=\sum_{\alpha\in\Delta_i}c_\alpha X_\alpha+c_{-\theta_i}\,t\,X_{-\theta_i}\in\fg\otimes\CC[t,t^{-1}].
\]
Every element of $\tOm_\cP$ normalizes $I$ and $T$, and it permutes the root vectors of $\cP$. It also preserves each connected component of the affine Dynkin diagram. Hence it fixes each $E_i$.

Let $H$ be the least common multiple of the $h_i$, so that $H\cdot\rho^\vee/h$ lies in $\X_*(T^{ad})$, and let $\gamma_H$ be the generator of $\Gal(\cK_H/\cK)$ from Section~\ref{subsec:galois}. In the notation of that section,
\[
t^{\rho^\vee/h}=\prod_i\rho_i^\vee(t^{1/h_i})\in T^{ad}(\cK_H),
\]
and by \eqref{eq:monodromy} its monodromy is the principal element:
\[
\gamma_H(t^{\rho^\vee/h})=x_0\,t^{\rho^\vee/h}.
\]
Conjugation by $t^{\rho^\vee/h}$ multiplies a root vector $X_\beta$ of the $i$-th factor by $t^{\langle\beta,\rho_i^\vee\rangle/h_i}$. Since $\langle\alpha,\rho_i^\vee\rangle=1$ for $\alpha\in\Delta_i$ and $\langle\theta_i,\rho_i^\vee\rangle=h_i-1$, every term of $E_i$ gets multiplied by $t^{1/h_i}$; for the last term, $t\cdot t^{-(h_i-1)/h_i}=t^{1/h_i}$. In other words,
\begin{equation}\label{eq:E}
\Ad(t^{\rho^\vee/h})\,E_i=t^{1/h_i}\,E_{0,i}\qquad\text{for each }i.
\end{equation}

\begin{prop}\label{prop:coxeter-torus}
Evaluation at $t=1$ gives an isomorphism
\[
\tOm_\cP\cong T_0^{\,c}.
\]
More precisely, every $g\in \tOm_\cP$ can be written uniquely as $g=\Ad(t^{-\rho^\vee/h})(\tau)$ with $\tau\in T_0^{\,c}$, and then $g|_{t=1}=\tau$.
\end{prop}

\begin{proof}
Let $g\in \tOm_\cP$, and put $\tau=\Ad(t^{\rho^\vee/h})(g)\in G(\cK_H)$.

\emph{Step 1: $\tau$ is a constant element of $T_0$.} Since $g$ fixes $E_i$, equation \eqref{eq:E} shows that $\tau$ fixes $E_{0,i}$ for every $i$, and hence fixes $E_0$. By Fact~\ref{fact:kostant}(1), the centralizer of $E_0$ in $G$ is the torus $T_0$, so $\tau\in T_0(\cK_H)$. Since $\tOm_\cP$ is finite, $\tau$ has finite order. The points of finite order of $T_0(\cK_H)$ are constant, so $\tau\in T_0(\CC)$.

\emph{Step 2: $\tau$ is fixed by $c$.} The element $g$ has coefficients in $\cK$, so it is fixed by $\gamma_H$. We compute $\gamma_H(g)$ in a second way. Since $\tau$ is constant, $\gamma_H$ only acts on the factors $t^{\mp\rho^\vee/h}$, which it multiplies by $x_0^{\mp1}$. Hence
\[
\gamma_H(g)=\Ad(t^{-\rho^\vee/h})\big(\Ad(x_0^{-1})\,\tau\big).
\]
Comparing with $g=\Ad(t^{-\rho^\vee/h})(\tau)$, we get $\Ad(x_0^{-1})\,\tau=\tau$. By the definition of $c$, this means that $\tau\in T_0^{\,c}$.

\emph{Step 3: $g\mapsto\tau$ is evaluation at $t=1$.} Setting $t^{1/H}=1$ in $g=\Ad(t^{-\rho^\vee/h})(\tau)$ gives $g|_{t=1}=\tau$. So $g\mapsto\tau$ is an injective homomorphism $\tOm_\cP\ra T_0^{\,c}$.

\emph{Step 4: both groups have the same order.} By Proposition~\ref{prop:Z-extension} and Corollary~\ref{cor:surjective}, the group $\tOm_\cP$ is an extension of $\Om_G$ by $Z(G)$. So
\[
|\tOm_\cP|=|Z(G)|\cdot|\Om_G|=|P^\vee/\X_*|\cdot|\X_*/Q^\vee|=|P^\vee/Q^\vee|.
\]
On the other hand, $c-1$ is an isogeny of $T_0$ with kernel $T_0^{\,c}$. Its degree is $|\det(1-c)|$, so by Fact~\ref{fact:det}
\[
|T_0^{\,c}|=|\det(1-c)|=|P^\vee/Q^\vee|.
\]
Hence the injective map of Step~3 is an isomorphism.
\end{proof}

\begin{rem}
In particular, $\tOm_\cP$ is commutative when $k=\CC$, since it embeds into a torus; Theorem~\ref{thm:duality}(1) extends this to positive characteristic. Moreover, the Tits lifts $n_{w_1}$ and $n_{w_\ell}$ in type $D_\ell$ with $\ell$ even commute as soon as they preserve the same affine pinning (compare Proposition~\ref{prop:generators} and Lemma~\ref{lem:pinned-D}).
\end{rem}

\subsection{The extension in lattice terms}
Fixed points of a Coxeter element have a simple description in terms of lattices. Write $\X_*=\X_*(T_0)$, and let $T_0^{sc}\subset G^{sc}$ be the preimage of $T_0$, so that $\X_*(T_0^{sc})=Q^\vee$. Both lattices span the same vector space $V=\X_*\otimes\CC$. We use two exponential maps,
\[
\exp\colon V\ra T_0^{sc}(\CC)\qquad\text{and}\qquad \pi\circ\exp\colon V\ra T_0(\CC),
\]
with kernels $Q^\vee$ and $\X_*$. By Lemma~\ref{lem:exp}(2), applied to $T_0$, the first map sends $\X_*$ onto $K(\CC)$. The action of $c$ on $V$ is compatible with both maps.

An element $\pi(\exp\xi)$ of $T_0$ is fixed by $c$ if and only if $\pi(\exp(c\,\xi))=\pi(\exp\xi)$, that is, if and only if $(c-1)\xi\in\X_*$. In $T_0^{sc}$, the elements $\exp(c\,\xi)$ and $\exp(\xi)$ then differ by the central element $\exp((c-1)\xi)\in K(\CC)$, which is trivial only if $(c-1)\xi\in Q^\vee$. We define
\[
\Lambda\colon T_0^{\,c}\ra\X_*/(1-c)\X_*,\qquad \pi(\exp\xi)\mapsto(c-1)\xi.
\]
This is well defined, because $\xi$ is determined modulo $\X_*$. It is injective, because $(c-1)\xi\in(c-1)\X_*$ forces $\xi\in\X_*$. It is surjective, because $c-1$ is invertible on $V$. So $\Lambda$ is an isomorphism.

\begin{prop}\label{prop:lattice}
Under the isomorphisms $\tOm_\cP\cong T_0^{\,c}\cong\X_*/(1-c)\X_*$, the exact sequence $1\ra Z(G)\ra \tOm_\cP\ra\Om_G\ra1$ becomes
\[
0\ra Q^\vee/(1-c)\X_*\ra\X_*/(1-c)\X_*\ra\X_*/Q^\vee\ra0,
\]
where $Q^\vee\subset\X_*$ is the coroot lattice and the second map is the projection.
\end{prop}

\begin{proof}
It suffices to show that the Kottwitz homomorphism $\tOm_\cP\ra\Om_G$ becomes the projection $\X_*/(1-c)\X_*\ra\X_*/Q^\vee$. Then the kernels correspond as well, so $Z(G)$ corresponds to $Q^\vee/(1-c)\X_*$.

We compute the Kottwitz homomorphism by monodromy, as in Lemma~\ref{lem:kummer}, using lifts over $\cK_H$. This is allowed because $H$ is a multiple of the exponent of $\Om$: for each simple factor, the exponent of $P^\vee/Q^\vee$ divides the Coxeter number, as one sees from the tables in \cite{Bou02}.

Let $g\in \tOm_\cP$. By Proposition~\ref{prop:coxeter-torus}, $g=\Ad(t^{-\rho^\vee/h})(\tau)$ with $\tau=g|_{t=1}\in T_0^{\,c}$. Write $\tau=\pi(\exp\xi)$ with $\xi\in V$, so that $(c-1)\xi\in\X_*$ and $\Lambda(\tau)$ is the class of $(c-1)\xi$. Then
\[
\tilde g=\Ad(t^{-\rho^\vee/h})(\exp\xi)\in G^{sc}(\cK_H)
\]
is a lift of $g$. Put $z=\exp\big((c-1)\xi\big)$. This element lies in $K(\CC)$, because $(c-1)\xi\in\X_*$, so it is central. The computation in Step~2 of the proof of Proposition~\ref{prop:coxeter-torus} gives $\gamma_H(\tilde g)=\Ad(t^{-\rho^\vee/h})\big(\Ad(x_0^{-1})\exp\xi\big)$. Here $\Ad(x_0^{-1})$ acts on $T_0^{sc}$ through $c$, so $\Ad(x_0^{-1})\exp\xi=\exp(c\,\xi)=\exp(\xi)\,z$. Since $z$ is central, we get
\[
\gamma_H(\tilde g)=\Ad(t^{-\rho^\vee/h})\big(\exp(\xi)\,z\big)=\tilde g\,z.
\] So the monodromy of $\tilde g$ is $z$. Under $K(\CC)\cong\X_*/Q^\vee$, the element $z$ corresponds to the class of $(c-1)\xi$, and by Lemma~\ref{lem:kummer} this class is $\kappa_G(g)$. Hence $\kappa_G(g)$ is the image of $\Lambda(\tau)$ under the projection $\X_*/(1-c)\X_*\ra\X_*/Q^\vee$, as claimed.
\end{proof}

Conjugating $T_0$ to $T$ by an element of $G$ identifies $\X_*(T_0)$ with $\X_*(T)$, compatibly with the coroot lattices, and turns $c$ into a Coxeter element of $W$. Since all Coxeter elements of $W$ are conjugate, we may take it to be any given one. So Proposition~\ref{prop:lattice} proves Theorem~\ref{thm:duality}(1) for $k=\CC$.

\subsection{The pairing}\label{subsec:pairing}
Apply Proposition~\ref{prop:lattice} to $G^\vee$ as well. Since all Coxeter elements are conjugate, we may use the same $c$, acting on $\X_*(T^\vee)=\X^*(T)$. This gives
\[
\tOm_\cP\cong\X_*/(1-c)\X_*,\qquad \tOm_{\cP^\vee}\cong\X^*/(1-c)\X^*,
\]
compatibly with the two exact sequences.

\begin{prop}\label{prop:linking}
The pairing
\[
\X_*/(1-c)\X_*\ \times\ \X^*/(1-c)\X^*\ra\CC^\times,\qquad (x,y)\mapsto\exp\big(2\pi i\,\langle(c-1)^{-1}x,\,y\rangle\big),
\]
is well defined and perfect.
\end{prop}

\begin{proof}
Put $A=c-1$, which is invertible on $\X_*\otimes\Q$. Since $\langle cx,cy\rangle=\langle x,y\rangle$, we have $\langle Ax,y\rangle=\langle x,(c^{-1}-1)y\rangle$, and $(c^{-1}-1)\X^*=A\X^*$ because $c^{-1}-1=-c^{-1}A$. If $x=Ax'$ with $x'\in\X_*$, then $\langle A^{-1}x,y\rangle=\langle x',y\rangle$ is an integer. If $y=(c^{-1}-1)y'$ with $y'\in\X^*$, then $\langle A^{-1}x,y\rangle=\langle x,y'\rangle$ is an integer. Hence the pairing is well defined. Both groups have order $|\det A|$, so it suffices that the left kernel is trivial. If $\langle A^{-1}x,y\rangle\in\Z$ for all $y\in\X^*$, then $A^{-1}x\in\X_*$, since $\X_*$ and $\X^*$ are dual lattices. So $x\in A\X_*$.
\end{proof}

Explicitly, if $g\in \tOm_\cP$ corresponds to $\pi(\exp\xi)\in T_0^{\,c}$ and $g^\vee\in \tOm_{\cP^\vee}$ corresponds to the class of $y\in\X^*$, then
\[
\langle g,g^\vee\rangle=\exp\big(2\pi i\langle\xi,y\rangle\big)=y\big(g|_{t=1}\big).
\]

\begin{proof}[Proof of Theorem~\ref{thm:duality}(2) for $k=\CC$]
We use the pairing of Proposition~\ref{prop:linking}. By Proposition~\ref{prop:lattice}, $Z(G)$ corresponds to $Q^\vee/(1-c)\X_*$. Its annihilator under the pairing of Proposition~\ref{prop:linking} consists of the classes of $y\in\X^*$ with
\[
\langle(c-1)^{-1}x,y\rangle\in\Z\qquad\text{for all }x\in Q^\vee.
\]
By Fact~\ref{fact:det}, $(c-1)^{-1}Q^\vee=P^\vee$, so the condition says $y\in(P^\vee)^*=Q$. So the annihilator is $Q/(1-c)\X^*$, which corresponds to $Z(G^\vee)$ by Proposition~\ref{prop:lattice} for $G^\vee$.

To identify the induced pairings, let $z=\pi(\exp\xi)\in Z(G)$ with $\xi\in P^\vee$, so $x=(c-1)\xi$. For $y\in\X^*$ we get
\[
\exp\big(2\pi i\langle(c-1)^{-1}x,y\rangle\big)=\exp\big(2\pi i\langle\xi,y\rangle\big)=y(z).
\]
This is the standard pairing $Z(G)\times\X^*/Q\ra\CC^\times$. The other pairing is symmetric. Finally, an extension of finite abelian groups splits if and only if its Pontryagin dual splits.
\end{proof}

\subsection{Fields of positive characteristic}\label{subsec:char-p}
We now remove the assumption $k=\CC$. The argument of Section~\ref{subsec:coxeter-torus} needs the sum of the root vectors of an affine pinning to be regular semisimple. This holds when $\chara k$ exceeds the Coxeter number of every simple factor \cite[Remark~3.20]{DLYZ}, but it can fail otherwise. Instead, we lift pinned elements to characteristic $0$. This only needs $p\nmid|Z(G)|$, so that the center does not change under reduction.

\begin{prop}\label{prop:reduction}
Let $k$ be an algebraically closed field of characteristic $p>0$ with $p\nmid|Z(G)|$, and let $\cP$ be an affine pinning of $LG_k$. Then there are an algebraically closed field $L$ of characteristic $0$ and an affine pinning $\cP_L$ of $LG_L$ such that the extensions
\[
1\ra Z(G)(L)\ra\tOm_{\cP_L}\ra\Om\ra1\qquad\text{and}\qquad1\ra Z(G)(k)\ra\tOm_{\cP}\ra\Om\ra1
\]
are isomorphic. Moreover, the isomorphism can be chosen so that it restricts to the reduction map $Z(G)(L)\ra Z(G)(k)$ on the kernels.
\end{prop}

\begin{proof}
Let $\cO$ be the ring of Witt vectors of $k$, let $L$ be an algebraic closure of its fraction field, and let $\cO_L$ be the integral closure of $\cO$ in $L$. The residue field of $\cO_L$ is $k$.

By Lemma~\ref{lem:conj-pinning}, we may assume that $\cP=\cP_c$ for some $c\in(k^\times)^{\Daff}$, where the reference pinning is the pinning of Corollary~\ref{cor:surjective}(1). This reference pinning is defined over $\Z$. Lift $c$ to $\tilde c\in(\cO^\times)^{\Daff}$, and let $\cP_L$ be the affine pinning $\cP_{\tilde c}$ of $LG_L$.

By Lemma~\ref{lem:conj-pinning} and the proof of Corollary~\ref{cor:surjective}, every element of $\tOm_{\cP_L}$ has the form $g_\eta\tau$. Here $g_\eta\in N_G(T)(\Z[t,t^{-1}])$, and $\tau\in T(L)$ satisfies $a(\tau)=\tilde c_{\eta(a)}/\tilde c_a$ for all $a\in\Daff$. These values are units in $\cO$. So the image of $\tau$ in $T^{ad}$ is integral, and hence $\tau\in T(\cO_L)$, because $\cO_L$ is integrally closed in $L$. Thus $\tOm_{\cP_L}$ lies in $N_G(T)(\cO_L[t,t^{-1}])$.

Reduction modulo the maximal ideal of $\cO_L$ now gives a homomorphism $\rho\colon\tOm_{\cP_L}\ra\tOm_{\cP}$, because $\cP_L$ reduces to $\cP$. It is compatible with $\kappa_G$, and on $Z(G)(L)=Z(G)(\cO_L)$ it is the reduction map to $Z(G)(k)$, which is an isomorphism because $Z(G)$ is finite étale over $\cO_L$. So $\rho$ is injective, since its kernel lies in $\ker\kappa_G\cap\tOm_{\cP_L}=Z(G)(L)$. By Proposition~\ref{prop:Z-extension} and Corollary~\ref{cor:surjective}, both groups have order $|Z(G)|\cdot|\Om|$, so $\rho$ is an isomorphism.
\end{proof}

\begin{proof}[Proof of Theorem~\ref{thm:duality}]
If $\chara k=0$, the proofs of Sections~\ref{subsec:kostant}--\ref{subsec:pairing} work verbatim over $k$. Suppose that $\chara k=p>0$.

For (1), assume $p\nmid|Z(G)|$. Proposition~\ref{prop:reduction} identifies the extension over $k$ with an extension over an algebraically closed field of characteristic $0$, and (1) follows from the case of characteristic $0$.

For (2), assume $p\nmid|Z(G)|\cdot|\pi_1(G)|$. Since $|Z(G^\vee)|=|\pi_1(G)|$, Proposition~\ref{prop:reduction} applies to both $G$ and $G^\vee$. It identifies both extensions with extensions in characteristic $0$, and on the centers these identifications are reduction of roots of unity. The standard pairings are compatible with reduction of roots of unity. So (2) follows from the case of characteristic $0$.
\end{proof}

\subsection{A criterion for pinnability}
Theorem~\ref{thm:dual}(1) turns pinnability into a question about lattices.

\begin{cor}\label{cor:criterion}
Let $k$ be an algebraically closed field with $\chara k\nmid|Z(G)|$. Then $G$ is pinnable over $k$ if and only if the sequence
\[
0\ra Q^\vee/(1-c)\X_*\ra\X_*/(1-c)\X_*\ra\X_*/Q^\vee\ra0
\]
splits. This holds if and only if $\X_*/(1-c)\X_*\cong Z(G)\times\pi_1(G)$ as abstract groups.
\end{cor}

\begin{proof}
The first statement follows from Theorem~\ref{thm:dual}(1). For the second, let $0\ra A\ra M\ra B\ra0$ be an exact sequence of finite abelian groups with $M\cong A\oplus B$. Applying $\Hom(B,-)$ gives an exact sequence
\[
0\ra\Hom(B,A)\ra\Hom(B,M)\ra\Hom(B,B).
\]
Since $|\Hom(B,M)|=|\Hom(B,A)|\cdot|\Hom(B,B)|$, the last map is surjective. So $\id_B$ lifts to a homomorphism $B\ra M$, and the sequence splits. The converse is clear.
\end{proof}

The same computation as in the next example, for $(SL_2\times Sp_{2n})/\mu_2$, gives $\Z/4\Z$ when $n$ is even and $(\Z/2\Z)^2$ when $n$ is odd. This recovers Theorem~\ref{thm:B}(2) over algebraically closed fields of characteristic $\neq2$, without Tits lifts. Here is a group that is not covered by Theorem~\ref{thm:B}.

\begin{exmp}\label{exmp:sp4sp6}
Let $G=(Sp_4\times Sp_6)/\mu_2$, with $\mu_2$ embedded diagonally. Then $Z(G)\cong\Z/2\Z$ and $\pi_1(G)\cong\Z/2\Z$. We show that $\X_*/(1-c)\X_*\cong\Z/4\Z$, so $G$ is not pinnable.

Use coordinates $e_1,e_2$ for $Sp_4$ and $f_1,f_2,f_3$ for $Sp_6$. Write elements of $\X_*\otimes\Q$ as
\[
(b;d)=b_1e_1+b_2e_2+d_1f_1+d_2f_2+d_3f_3.
\]
Then $Q^\vee=\Z^2\oplus\Z^3$, and $\X_*$ consists of the $(b;d)$ such that either
\[
b\in\Z^2\text{ and }d\in\Z^3,\qquad\text{or}\qquad b\in(\tfrac12,\tfrac12)+\Z^2\text{ and }d\in(\tfrac12,\tfrac12,\tfrac12)+\Z^3.
\]
Take the Coxeter elements $e_1\mapsto e_2\mapsto-e_1$ of $C_2$ and $f_1\mapsto f_2\mapsto f_3\mapsto-f_1$ of $C_3$. Then
\[
(1-c)(b;d)=(b_1+b_2,\ b_2-b_1;\ d_1+d_3,\ d_2-d_1,\ d_3-d_2).
\]
The element $v=(\frac12,\frac12;\frac12,\frac12,\frac12)$ generates $\pi_1(G)$. Suppose that $2v=(1,1;1,1,1)$ equals $(1-c)(b;d)$ for some $(b;d)\in\X_*$.
\begin{itemize}
\item The equations $b_1+b_2=1$ and $b_2-b_1=1$ give $b=(0,1)$, which is integral.
\item The equations $d_1+d_3=1$, $d_2-d_1=1$ and $d_3-d_2=1$ give $d=(-\frac12,\frac12,\frac32)$, which is not integral.
\end{itemize}
So $(b;d)\notin\X_*$, a contradiction. Hence $v$ has order $4$ in $\X_*/(1-c)\X_*$. This group has order $|Z(G)|\cdot|\pi_1(G)|=4$, so it is cyclic of order $4$, and the extension does not split.

This agrees with a direct computation. The pinned lifts of Lemma~\ref{lem:pinned-tits} in the two factors have squares $\Id$ and $-\Id$ (proof of Proposition~\ref{prop:generators}). So the pinned lift of the generator of $\Om_G$ has square $(\Id,-\Id)$, which is not in $\mu_2$. By Theorem~\ref{thm:dual}, the dual group $(Spin_5\times Spin_7)/\mu_2$ is not pinnable either.
\end{exmp}

\section{Centralizers of semisimple elements}\label{sec:centralizers}

Let $s$ be a semisimple element of $G$. Its centralizer $G_s$ is often disconnected, and its component group plays a role in endoscopy. In this section we prove Theorem~\ref{thm:cent}: the section of Theorem~\ref{thm:A}, restricted to a suitable subgroup of $\Om$ and evaluated at $t=1$, splits the component group of $G_s$.

By Lemma~\ref{lem:finite-order}, we may first replace $s$ by an element of finite order with the same centralizer. The rest of the proof uses the geometry of the fundamental alcove. We write $s=\exp(\chi)$ with $\chi$ in the closed alcove. The walls of the alcove through $\chi$ give a basis of the root system of $G_s^{\circ}$ (Lemma~\ref{lem:walls}), and the stabilizer of $\chi$ in $\Om$ is the component group of $G_s$ (Proposition~\ref{prop:pi0}). The section of Theorem~\ref{thm:A} permutes these walls, so after evaluation at $t=1$ it preserves a Borel subgroup of $G_s^{\circ}$.

\subsection{Reduction to elements of finite order}
In this section, $k$ is algebraically closed of arbitrary characteristic; only Proposition~\ref{prop:centralizer} assumes $\chara k\neq2$. Let $s\in T(k)$. By a theorem of Steinberg \cite[II,~4.1]{SS70}, $G_s^{\circ}$ is generated by $T$ and the root subgroups $U_\alpha$ with $\alpha(s)=1$, and $G_s$ is generated by $G_s^{\circ}$ and the elements of $N_G(T)$ that centralize $s$. So the group $G_s$ is determined by the two sets
\[
\Phi_s=\{\alpha\in\Phi:\alpha(s)=1\},\qquad W_s=\{w\in W:w(s)=s\}.
\]

\begin{lem}\label{lem:finite-order}
Let $s\in T(k)$. There is an element $s'\in T(k)$ of finite order with $\Phi_{s'}=\Phi_s$ and $W_{s'}=W_s$. In particular, $G_{s'}=G_s$.
\end{lem}

\begin{proof}
Let $H\subset T$ be the Zariski closure of the subgroup generated by $s$. It is a diagonalizable group, so $H=H^{\circ}\times F$ for a finite group $F$ of order prime to $\chara k$. The coset $sH^{\circ}$ is irreducible, and its points of finite order are Zariski dense in it.

Every root $\alpha$ with $\alpha(s)=1$ is trivial on $H$. Every $w\in W$ with $w(s)=s$ fixes $H$ pointwise. So $\Phi_s\subset\Phi_{h}$ and $W_s\subset W_{h}$ for every $h\in H$.

It remains to choose $s'$ so that no further roots are trivial on $s'$, and no further elements of $W$ fix $s'$. For this, consider the finitely many closed subsets
\[
\{h\in sH^{\circ}:\alpha(h)=1\}\ (\alpha\notin\Phi_s),\qquad\{h\in sH^{\circ}:w(h)=h\}\ (w\notin W_s).
\]
None of them contains $s$, so they are proper closed subsets of the irreducible variety $sH^{\circ}$. Hence there is a point $s'\in sH^{\circ}$ of finite order outside all of them. Then $\Phi_{s'}=\Phi_s$ and $W_{s'}=W_s$.
\end{proof}

\subsection{Setup} Let $p=\chara k$, and let $\Q_{p'}\subset\Q$ be the subring of fractions with denominator prime to $p$ (with $\Q_{p'}=\Q$ if $p=0$). The roots of unity of Section~\ref{subsec:galois} give an isomorphism between $\Q_{p'}/\Z$ and the group of roots of unity of order prime to $p$ in $k^\times$. They give a homomorphism
\[
\exp\colon\X_*(T)\otimes\Q_{p'}\ra T(k)
\]
with kernel $\X_*(T)$, whose image is the torsion subgroup of $T(k)$. For $k=\CC$ this is the map of Section~\ref{subsec:groups}.

Let $s\in G(k)$ be semisimple. We may assume $s\in T(k)$, and by Lemma~\ref{lem:finite-order} we may assume that $s$ has finite order. This order is prime to $p$, since $k^\times$ has no element of order $p$. Choose $\chi\in\X_*(T)\otimes\Q_{p'}$ with $\exp(\chi)=s$. For $g\in\extw$ we have $\exp(g\cdot\chi)=\pr(g)(s)$. The closure $\bar A$ of $A$ meets every $\Wa$-orbit in $V$ exactly once \cite[Chapter~V, \S3.3, Theorem~2]{Bou02}. Replacing $s$ by a $W$-conjugate, we may assume $\chi\in\bar A$. We put
\[
\Phi_\chi=\{\alpha\in\Phi:\langle\alpha,\chi\rangle\in\Z\},\qquad \Delta_{\aff,\chi}=\{a\in\Daff:a(\chi)=0\},
\]
and we let $\Om_\chi$ be the stabilizer of $\chi$ in $\Om$. As in Section~\ref{subsec:W}, for an affine root $a=\alpha+m\delta$ we write $\pi(a)=\alpha$ for its gradient. Since $\alpha(s)=1$ if and only if $\langle\alpha,\chi\rangle\in\Z$, we have $\Phi_\chi=\Phi_s$, the set of roots of $G_s^{\circ}$ with respect to $T$.

\subsection{Walls and components}
The next lemma is the slogan ``every pseudo-Levi subgroup of $G$ is the Levi quotient of a standard parahoric subgroup of the loop group''.

\begin{lem}\label{lem:walls}
Let $\chi\in\bar A$.
\begin{enumerate}[(1)]
\item The stabilizer of $\chi$ in $\Wa$ is generated by the reflections $s_a$ with $a\in\Delta_{\aff,\chi}$.
\item The stabilizer of $\chi$ in $\extw$ is the semidirect product of its stabilizer in $\Wa$ and $\Om_\chi$.
\item The set $\pi(\Delta_{\aff,\chi})$ is a basis of the root system $\Phi_\chi$.
\end{enumerate}
\end{lem}

\begin{proof}
Part (1) is \cite[Chapter~V, \S3.3, Proposition~1]{Bou02}. For (2), let $g=v\om$ fix $\chi$, with $v\in\Wa$ and $\om\in\Om$. Then $\om(\chi)\in\bar A$ and $v(\om(\chi))=\chi$. Since $\bar A$ meets each $\Wa$-orbit once, $\om(\chi)=\chi$.

For (3), let $\Phi_{\aff,\chi}$ be the set of affine roots vanishing at $\chi$. The map $a\mapsto\pi(a)$ is a bijection $\Phi_{\aff,\chi}\ra\Phi_\chi$, with inverse $\alpha\mapsto\alpha-\langle\alpha,\chi\rangle\delta$.

Now let $a\in\Phi_{\aff,\chi}$ be positive on $A$, and write $a=\sum_{b\in\Daff}n_b\,b$ with integers $n_b\geq0$. Evaluating at $\chi$ gives $0=\sum_bn_b\,b(\chi)$. Since $b(\chi)\geq0$ for every $b\in\Daff$, we get $n_b=0$ whenever $b(\chi)>0$. So $a$ is a non-negative combination of $\Delta_{\aff,\chi}$. Applying $\pi$, we see that every root of $\Phi_\chi$ is a non-negative or a non-positive integral combination of $\pi(\Delta_{\aff,\chi})$.

It remains to show that $\pi(\Delta_{\aff,\chi})$ is linearly independent. For each simple factor, the simple affine roots $b$ satisfy $\sum_bn_b\,b=\delta$ with all $n_b>0$, and $\delta$ is the constant function $1$. Evaluating at $\chi$ shows that $\Delta_{\aff,\chi}$ does not contain all simple affine roots of any simple factor. Let $\alpha_1,\dots,\alpha_\ell$ be the simple roots of a simple factor and $\theta=\sum_in_i\alpha_i$ its highest root. The gradients of the simple affine roots of this factor are $\alpha_1,\dots,\alpha_\ell$ and $-\theta$. Since every $n_i$ is positive, any $\ell$ of these $\ell+1$ vectors are linearly independent. Hence $\pi(\Delta_{\aff,\chi})$ is a basis of $\Phi_\chi$.
\end{proof}

\begin{prop}\label{prop:pi0}
The map $\pr$ induces an isomorphism
\[
\Om_\chi\cong\pi_0(G_s).
\]
Explicitly, $\om\in\Om_\chi$ corresponds to the component of $G_s$ that contains any lift of $\pr(\om)$.
\end{prop}

\begin{proof}
An element $w\in W$ fixes $s$ if and only if $w\chi-\chi\in\X_*(T)$. This holds if and only if $w\,t^{\mu}$ fixes $\chi$ for some $\mu$, and then $\mu$ is unique. So $\pr$ is an isomorphism from the stabilizer of $\chi$ in $\extw$ onto $W_s$.
By Lemma~\ref{lem:walls}(1) and (3), it maps the stabilizer of $\chi$ in $\Wa$ onto the Weyl group of $\Phi_\chi$, which is the Weyl group of $G_s^{\circ}$. All maximal tori of $G_s^{\circ}$ are conjugate under $G_s^{\circ}$, so $G_s=G_s^{\circ}N_{G_s}(T)$. Hence $\pi_0(G_s)\cong W_s/W(G_s^{\circ})$. By Lemma~\ref{lem:walls}(2), this quotient is isomorphic to $\Om_\chi$.
\end{proof}

\subsection{The finite part map}
The splitting in Theorem~\ref{thm:cent} is obtained from a section of $\kappa_G$ by forgetting its translation part. Since $T(\cK)=t^{\X_*(T)}\times T(\cO)$, every element of $N_G(T)(\cK)$ can be written uniquely as $n_w\,t^{\mu}\,\tau$ with $w\in W$, $\mu\in\X_*(T)$ and $\tau\in T(\cO)$. We put
\[
\varphi(n_w\,t^{\mu}\,\tau)=n_w\,\tau(0)\in N_G(T)(k).
\]
In words, $\varphi$ forgets $t^\mu$ and replaces $\tau$ by its constant term.

For the section $\psi$ of Theorem~\ref{thm:A}, the torus part is a constant element of $T(\Z)$, so $\varphi(\psi(\om))$ is simply $\psi(\om)$ evaluated at $t=1$.

\begin{lem}\label{lem:finite-part}
The map $\varphi\colon N_G(T)(\cK)\ra N_G(T)(k)$ is a group homomorphism. On the image of $N_G(T)(\Z[t,t^{-1}])$, it is given by $t\mapsto1$.
\end{lem}

\begin{proof}
Let $g_j=n_{w_j}t^{\mu_j}\tau_j$ for $j=1,2$. Moving $n_{w_2}$ to the left, we get
\[
g_1g_2=n_{w_1}n_{w_2}\,t^{w_2^{-1}\mu_1+\mu_2}\,w_2^{-1}(\tau_1)\,\tau_2.
\]
Write $n_{w_1}n_{w_2}=n_{w_1w_2}\,\epsilon$ with $\epsilon\in T(\Z)$ (Lemma~\ref{lem:adrian}(1)). Then
\[
\varphi(g_1g_2)=n_{w_1}n_{w_2}\,w_2^{-1}(\tau_1(0))\,\tau_2(0)=n_{w_1}\tau_1(0)\,n_{w_2}\tau_2(0)=\varphi(g_1)\varphi(g_2).
\]
For the second statement, write $n\in N_G(T)(\Z[t,t^{-1}])$ as $n=\tau\,n_w\,t^{\mu}=n_w\,t^{\mu}\,w^{-1}(\tau)$ with $\tau\in T(\Z)$. Then $\varphi(n)=n_w\,w^{-1}(\tau)=\tau\,n_w$, which is the image of $n$ under $t\mapsto1$.
\end{proof}

\subsection{Proof of Theorem~\ref{thm:cent}}
We can now state and prove the main result of this section.

\begin{thm}\label{thm:cent-body}
Let $\sigma\colon\Om\ra LG(k)$ be a homomorphic section of $\kappa_G$ whose image lies in $N_G(T)(\cK)$. Then the restriction of $\varphi\circ\sigma$ to $\Om_\chi$ is a splitting of
\[
1\ra G_s^{\circ}\ra G_s\ra\pi_0(G_s)\ra1.
\]
This splitting preserves $T$ and the Borel subgroup $B_s\subset G_s^{\circ}$ whose simple roots are $\pi(\Delta_{\aff,\chi})$.

If, moreover, $\sigma(\Om)$ preserves an affine pinning $(I,T,(x_a))$, then the splitting preserves the pinning
\[
\big(B_s,\ T,\ (\bar x_a)_{a\in\Delta_{\aff,\chi}}\big)
\]
of $G_s^{\circ}$. Here $\bar x_a=u_\alpha(c)$ if $x_a=u_\alpha(c\,t^m)$ and $a=\alpha+m\delta$.
\end{thm}

\begin{proof}
The map $\varphi$ is a homomorphism by Lemma~\ref{lem:finite-part}. For $\om\in\Om_\chi$, the element $\varphi(\sigma(\om))$ lifts $\pr(\om)$, which fixes $s$. So it lies in $N_{G_s}(T)$, and by Proposition~\ref{prop:pi0} it lies in the component corresponding to $\om$. This gives the splitting.

Next, $\sigma(\om)$ normalizes $I$, so $\sigma(\om)U_a\sigma(\om)^{-1}=U_{\om\cdot a}$ for all $a\in\Daff$. Write $\sigma(\om)=n_w\,t^{\mu}\,\tau$, and let $a=\alpha+m\delta$ and $c\in k$. Conjugating successively by $\tau$, $t^{\mu}$ and $n_w$, we get
\[
\sigma(\om)\,u_\alpha(c\,t^m)\,\sigma(\om)^{-1}=n_w\,u_\alpha\big(\alpha(\tau)\,c\,t^{m+\langle\alpha,\mu\rangle}\big)\,n_w^{-1}.
\]
In the same way,
\[
\varphi(\sigma(\om))\,u_\alpha(c)\,\varphi(\sigma(\om))^{-1}=n_w\,u_\alpha\big(\alpha(\tau(0))\,c\big)\,n_w^{-1}.
\]
The second formula shows that $\varphi(\sigma(\om))$ maps $U_{\pi(a)}$ onto $U_{\pi(\om\cdot a)}$. Since $\om\in\Om_\chi$ permutes $\Delta_{\aff,\chi}$, and $\pi(\Delta_{\aff,\chi})$ is the set of simple roots of $B_s$ by Lemma~\ref{lem:walls}(3), the element $\varphi(\sigma(\om))$ preserves $B_s$ and $T$.

Finally, suppose that $\sigma(\Om)$ preserves the affine pinning $(I,T,(x_a))$. Let $\om\in\Om_\chi$ and $a=\alpha+m\delta\in\Daff$, and write $x_a=u_\alpha(c\,t^m)$. Since $\sigma(\om)$ maps $x_a$ to $x_{\om\cdot a}$, the first formula shows that $\alpha(\tau)$ is a constant, so $\alpha(\tau)=\alpha(\tau(0))$. Comparing the two formulas, we see that $\varphi(\sigma(\om))$ maps $\bar x_a=u_\alpha(c)$ to $\bar x_{\om\cdot a}$. As $\om$ permutes $\Delta_{\aff,\chi}$, the splitting preserves the pinning $(B_s,T,(\bar x_a)_{a\in\Delta_{\aff,\chi}})$.
\end{proof}

Theorem~\ref{thm:cent}(1) follows by applying Theorem~\ref{thm:cent-body} to the section of Theorem~\ref{thm:A-body}. Theorem~\ref{thm:cent}(2) follows by applying it to a pinned splitting. Every element of $LG(k)$ that normalizes $I$ and $T$ lies in $N_G(T)(\cK)$, so the theorem applies.

\subsection{A counterexample}
Finally, we prove Theorem~\ref{thm:cent}(3). The groups of Theorem~\ref{thm:B}(2) are not pinnable, and this shows up in a centralizer.

\begin{prop}\label{prop:centralizer}
Let $\chara k\neq2$, let $n\geq2$ be even, and let $G=(SL_2\times H)/\langle(-\Id,z)\rangle$ be one of the groups in Proposition~\ref{prop:B2-body}. Let $i\in k$ with $i^2=-1$. Let $s\in G(k)$ be the image of $(\diag(i,-i),s_H)$, where $s_H\in T_H(k)$ is defined as follows:
\begin{itemize}
\item if $H=Sp_{2n}$, then $s_H=\diag(i,\dots,i,-i,\dots,-i)$ in the coordinates of Example~\ref{exmp:sp4};
\item if $H=Spin_{2n+1}$, then $s_H=L_1^\vee(i)$, where $L_1^\vee=2e_1$ is the coroot of the short root $L_1$.
\end{itemize}
Then $s$ has order $2$, and the group $\pi_0(G_s)$ has order $2$. No splitting of
\[
1\ra G_s^{\circ}\ra G_s\ra\pi_0(G_s)\ra1
\]
preserves a pinning of $G_s^{\circ}$.
\end{prop}

\begin{proof}
In both cases $s_H^2=z$. So $s^2$ is the image of $(-\Id,z)$, and $s$ has order $2$.

\emph{The component group.} Let $\sigma$ be the non-trivial Weyl element of $SL_2$, and let $w_H$ be the finite part of the generator of $\Om_{H^{ad}}$. So $w_H$ sends $e_j$ to $-e_{n+1-j}$ if $H=Sp_{2n}$, and $w_H=s_{L_1}$ if $H=Spin_{2n+1}$. In both cases $w_H(s_H)=zs_H$. So $w=(\sigma,w_H)$ sends $(\diag(i,-i),s_H)$ to $(\diag(-i,i),zs_H)$, and $w$ fixes $s$.

The roots that are trivial on $s$ are as follows. If $H=Sp_{2n}$, they are the roots $L_j-L_k$, which form a root system of type $A_{n-1}$ with simple roots $\alpha_1,\dots,\alpha_{n-1}$. The element $w$ swaps $\alpha_j$ and $\alpha_{n-j}$, and it fixes $\beta=\alpha_{n/2}$. If $H=Spin_{2n+1}$, they are the roots that do not involve $L_1$. They form a root system of type $B_{n-1}$ with simple roots $\alpha_2,\dots,\alpha_n$, all fixed by $w$. One checks that the stabilizer of $s$ in $W$ is $W(G_s^{\circ})\sqcup W(G_s^{\circ})\,w$. So $\pi_0(G_s)$ has order $2$, and $w$ is the only element of the non-trivial coset that preserves the simple roots above.

\emph{Strategy.} Suppose that $y\in G_s$ has order $2$, lies in the non-trivial component, and preserves a pinning of $G_s^{\circ}$. We derive a contradiction. Conjugating by $G_s^{\circ}$, we may assume that this pinning has torus $T$ and the simple roots listed above. Then $y\in N_G(T)$, and $y$ lifts $w$, because $w$ is the only element of the non-trivial coset that preserves these simple roots. Write $y$ as the image of $(y_1,y_2)$, where $y_1\in SL_2$ lifts $\sigma$ and $y_2=n_{w_H}\tau$ with $\tau\in T_H$. In both cases $n_{w_H}^2=1$ because $n$ is even (proof of Proposition~\ref{prop:generators}). Hence
\[
y_2^2=w_H(\tau)\,\tau.
\]
We get two conditions on $\tau$.
\begin{enumerate}[(a)]
\item Since $y_1^2=-\Id$, the condition $y^2=1$ says that $y_2^2=z$, that is, $w_H(\tau)\,\tau=z$.
\item Let $\alpha$ be a simple root of $G_s^{\circ}$ that is fixed by $w$. Then $y$ fixes the root vector of $\alpha$ in the pinning. Since $\alpha$ and $w_H(\alpha)=\alpha$ are simple roots of $H$, Lemma~\ref{lem:springer} shows that $n_{w_H}$ fixes $u_\alpha(c)$. So $y\,u_\alpha(c)\,y^{-1}=u_\alpha(\alpha(\tau)\,c)$, and the condition says that $\alpha(\tau)=1$.
\end{enumerate}
We now show that (a) and (b) are incompatible.

\emph{Case $H=Sp_{2n}$.} Write $\tau=\diag(t_1,\dots,t_n,t_n^{-1},\dots,t_1^{-1})$. The first $n$ diagonal entries of $w_H(\tau)\,\tau$ are $t_j/t_{n+1-j}$, so (a) says that $t_j=-t_{n+1-j}$ for all $j$. For the fixed simple root $\beta=\alpha_{n/2}$, this gives $\beta(\tau)=t_{n/2}/t_{n/2+1}=-1$. This contradicts (b).

\emph{Case $H=Spin_{2n+1}$.} Here $w_H=s_{L_1}$, and the simple roots $\alpha_2,\dots,\alpha_n$ of $G_s^{\circ}$ are fixed by $w$. By (b), $\alpha_j(\tau)=1$ for $2\leq j\leq n$. Since $\alpha_n=L_n$ and $\alpha_j=L_j-L_{j+1}$, this means that $L_j(\tau)=1$ for all $j\geq2$. Now evaluate the spin weight $\varpi_n=\frac12(L_1+\dots+L_n)$ on $w_H(\tau)\,\tau$:
\[
\varpi_n\big(w_H(\tau)\,\tau\big)=(\varpi_n+s_{L_1}\varpi_n)(\tau)=(L_2+\dots+L_n)(\tau)=1.
\]
But $\varpi_n(z)=-1$, so this contradicts (a).
\end{proof}

\bibliographystyle{amsalpha}
\bibliography{Kottwitz_paper}

\end{document}